\documentclass[11pt,reqno]{amsart}
\usepackage[T1]{fontenc}
\usepackage{lmodern}
\usepackage{mathrsfs}
\usepackage{amsmath,amsfonts,amssymb,amsthm,mathtools,booktabs}
\usepackage{microtype}
\usepackage[colorlinks=true,urlcolor=blue,citecolor=red,linkcolor=blue]{hyperref}
\newtheorem{theorem}{Theorem}[section]
\newtheorem{lemma}[theorem]{Lemma}
\newtheorem{proposition}[theorem]{Proposition}
\newtheorem{corollary}[theorem]{Corollary}
\newtheorem{conjecture}[theorem]{Conjecture}
\theoremstyle{definition}

\newtheorem{remark}[theorem]{Remark}
\numberwithin{equation}{section}
\newcommand{\R}{\mathbb R}
\newcommand{\Nm}{\operatorname{N}}
\newcommand{\Tr}{\operatorname{Tr}}
\newcommand{\lcm}{\operatorname{lcm}}
\newcommand{\Sq}{\operatorname{Sq}}
\newcommand{\Supp}{\operatorname{Supp}}
\newcommand{\ord}{\operatorname{ord}}
\newcommand{\num}{\operatorname{num}}
\newcommand{\hW}{h_{\mathrm W}}
\newcommand{\Q}{\mathbb Q}
\newcommand{\Z}{\mathbb Z}

\def\bg{\bigg}
\def\({\bg(}
\def\){\bg)}
\def\t{\text}

\def\mo{{\rm{mod}\ }}
\def\pmod#1{\ (\mo\ #1)}

\def\lcm{\rm{lcm}}
\def\ord{{\rm ord}}

\def\sm{\setminus}

\begin{document}
\hbox{Preprint, {\tt arXiv:2611.14605}}
\medskip

\title[$\mathbb Z$ is diophantine over $\mathbb Q$ with 26 unknowns]{$\mathbb Z$ is diophantine over $\mathbb Q$ with 26 unknowns}
\author[Zhi-Wei Sun]{Zhi-Wei Sun}
\address{School of Mathematics, Nanjing University, Nanjing 210093, People's Republic of China}
\email{zwsun@nju.edu.cn}

\keywords{Hilbert's tenth problem, existential definition of integers, rational numbers, Lucas sequences, elliptic formal groups, pole predicates.
\newline \indent 2020 {\it Mathematics Subject Classification}. Primary 11U05, 03D35; Secondary 11D99, 11G05.
\newline \indent
Supported by the National Natural Science Foundation of China (grant no. 12371004).}

\begin{abstract}
Recently, OpenAI proved that Hilbert's Tenth Problem over the field $\mathbb Q$
of rational numbers has a negative solution. Motivated by this work, G. Conant et al. showed that $\mathbb Z$ is diophantine over $\mathbb Q$, i.e., there is an integer-coefficient polynomial $F(x,y_1,\ldots,y_n)$
such that for any $x\in\mathbb Q$ we have 
$$ x\in\mathbb Z\quad\iff  \exists\,y_1,\ldots,y_{n}\in\mathbb Q\ [F(x,y_1,\ldots,y_{n})=0].$$
They constructed such a polynomial $F$ with $\deg F=4$ and $n=221$.  
In this paper, we prove the existence of such a polynomial $F$ with $n=26$.
We also show that $16$ existential unknowns suffice under the $abc$-conjecture.
\end{abstract}
\maketitle
\enlargethispage{2pt}

\section{Introduction}\label{sec:intro}

Hilbert's Tenth Problem over a ring asks for an algorithm deciding whether
an arbitrary polynomial over the ring has a zero in that ring. The negative
solution over $\mathbb Z$ was completed by Matiyasevich \cite{Mat} in 1970 following Davis,
Putnam and Robinson \cite{DPR}.

Let $\Q$ be the field of rational numbers.
Recently, OpenAI \cite{OAI} proved that Hilbert's Tenth Problem over $\Q$ also has a negative solution.
The proof combines elliptic curves, $p$-adic formal groups, Diophantine geometry, Selmer groups, height theory, and model-theoretic compactness. The author 
\cite{Sun86} proved further that there is no algorithm to decide for any $P(x_1,\ldots,x_{48})$
with integer coefficients whether $P(x_1,\ldots,x_{48})=0$
for some $x_1,\ldots,x_{48}\in\Q$. The proof combines an elliptic realization of integer
indices, a five-contact height estimate, simultaneous prime patterns, and
an economical algebraic compiler. 

A subset $A$ of $\Q$ is \emph{diophantine over $\mathbb Q$ with
$n$ unknowns} if there is a polynomial
$F\in\mathbb Z[X,Y_1,\ldots,Y_n]$ such that for any $x\in\Q$ we have
\begin{equation}\label{eq:def}
 x\in A\ \Longleftrightarrow\
 \exists (y_1,\ldots,y_n)\in\mathbb Q^n\ [F(x,y_1,\ldots,y_n)=0].
\end{equation}

In 1949  Robinson \cite{Rob} found a first-order definition of
$\mathbb Z$ in $\mathbb Q$, which was refined by Poonen \cite{Poonen} in 2009.
In 2016 Koenigsmann \cite{Koen} gave an existential definition of $\mathbb Q\setminus\mathbb Z$
in $\Q$.  Daans \cite{Daans} improved this by showing that $\Q\sm\Z$
is Diophantine over $\Q$ with $10$ unknowns; Sun \cite{Sun7} subsequently reduced $10$ to $7$.

Although OpenAI solved Hilbert's Tenth Problem over $\Q$, the paper \cite{OAI}
gives no single existential definition of $\Z$ in $\Q$. Quite recently,
 \cite{Glazer} obtained an
existential definition of $\mathbb Z$ in $\mathbb Q$. In particular,
In particular,
\cite[Theorem 1.2]{Glazer} gives a polynomial $F(x,y_1,\ldots,y_{221})$ with integer coefficients
of degree $4$ (which cannot be written explicitly with concrete coefficients) such that for any $x\in\Q$ we have
$$x\in\Z\iff \exists y_1,\ldots,y_{221}\in\Q\ [F(x,y_1,\ldots,y_{221})=0].$$

Burungale \cite{Burungale} subsequently gave an unconditional
Diophantine definition of $\mathbb Z$ in $\mathbb Q$ by a different
arithmetic route. His Theorem B realizes prescribed square classes on
congruent-number curves using the Tian--Yuan--Zhang analytic-rank criterion,
and his Theorem E proves a five-contact height estimate from quaternionic
surfaces and published modularity inputs. We use the latter estimate below,
after an explicit fractional-linear change of coordinate. His
Proposition 2.3 also gives the formal-logarithm comparison needed at the
large contact primes. We retain the two-structural-parameter predicate
$\Phi_0$ from \cite{Sun86} because it is already integrated into the
certificate and variable count below. Burungale provides an alternative
unconditional square-denominator mechanism, but we do not use that
alternative predicate in the present construction.

Our purpose is to decrease the number of existential unknowns in an
actual definition of $\mathbb Z$ significantly. The result is as follows.

\begin{theorem}\label{thm:main}
$\mathbb Z$ is
diophantine over $\mathbb Q$ with $26$ unknowns. In other words, 
there is a polynomial
\[
 F\in\mathbb Z[X,Y_1,\ldots,Y_{26}]
\]
such that for any $x\in\Q$ we have
$$x\in\Z\iff \exists (y_1,\ldots,y_{26})\in\Q^{26}\,[F(x,y_1,\ldots,y_{26})=0].$$ 
\end{theorem}

We also record a conditional refinement motivated by \cite{KPV}. We use the following standard
form of the $abc$-conjecture: for every $\varepsilon>0$ there is a constant
$C_\varepsilon>0$ such that, whenever $a,b,c$ are nonzero coprime integers
with $a+b=c$, one has
\[
 \max\{|a|,|b|,|c|\}\le C_\varepsilon
 \operatorname{rad}(abc)^{1+\varepsilon},
\]
where $\operatorname{rad}(m)$ is the product of all the distinct prime
divisors of an nonzero integer $m$.

\begin{theorem}\label{thm:abc16}
Under the $abc$-conjecture, $\mathbb Z$ is diophantine over
$\mathbb Q$ with $16$ unknowns.
\end{theorem}

Theorem~\ref{thm:abc16} is proved in Section~\ref{sec:abc}. Its new
input is only the $abc$-conjecture; all other arithmetic ingredients are
those already used for Theorem~\ref{thm:main}.

In our proof of Theorem \ref{thm:main}, we use the specific generator
certificate \cite[Proposition 3.1]{Glazer}, the odd-index certificate
\cite[Proposition J.7]{Glazer}, and the aggregate trace-height theorem
\cite[Section A38]{Glazer}. The rank-one and full rational two-torsion
properties of the same curve are also proved independently in
\cite[Proposition 2.1]{Burungale}. The five-contact height estimate
used here is obtained from \cite[Theorem E]{Burungale} by the explicit
fractional-linear transformation in Section~\ref{sec:inputs}; thus we
do not use the five-contact theorem of \cite{Glazer}. At large good
primes the formal-logarithm comparison is supplied by
\cite[Proposition 2.3]{Burungale}.
The standard fixed-multiplication formulas also appear in
\cite[Section 2.2]{Glazer}. From \cite{Sun86} we use
Proposition 2.3, Lemma 3.1, Theorem 4.1, Lemma 8.1, and the
finite-algebra compilation in Sections 10.1--10.6; its formal-logarithm
lemma may equivalently be replaced at the large contact primes by the
Burungale proposition just cited. Our contribution is the reduction
from these known facts to the $26$-unknown existential definition.

For the remainder of this introduction, $P$ denotes the fixed
infinite-order generator used in the construction; $Q$ is an auxiliary
cyclic point; $N,S,T,Z$ are the four cyclic root points; and $H$ is the
common comparison point.  The notation $[m]Q$ denotes scalar
multiplication on an elliptic curve, $x(Q)$ denotes the abscissa of an
affine point, and $\hW$ denotes the absolute logarithmic Weil height.
Section~\ref{sec:inputs} gives the full conventions before the proof
uses these symbols in a technical argument.

Here are the variable reductions.  Each of the nine required points of
the cyclic subgroup is represented by the two rational components of
the abscissa of a doubling preimage.  Its ordinate is a finite
algebraic choice.  For a positive integer $M$ fixed later, scalar multiplication computes
$R_{\rm sc}=[M]Q$, and the second global root is
$S=P+R_{\rm sc}$, so neither point costs a further pair of free
coordinates.  We impose finite algebraic odd-index certificates on the
other three global roots $N,T,Z$.  Hence all four root indices are odd.

For the four nonzero root indices, put
\[
 T_{\rm ag}=\frac12\Tr_{K/\Q}
 \bigl(x(N)+x(S)+x(T)+x(Z)\bigr)>0.
\]
By \cite[Section A38]{Glazer}, there are fixed constants $M_1>0$ and $M_2\ge0$ such that 
 $B^2\le M_1\hW(T_{\rm ag})+M_2$ with $B=\max\{|n|,|s|,|t|,|z|\}$.
Thus the two positive rational square scales
\[
 A_0=1\ \ \t{and}\ \ A_1=T_{\rm ag}^2
\]
control all four root heights.  Consequently, only two ordinary contact
predicates are needed, and they share one rational slope.  The scalar
height certificate is applied to $R_{\rm sc}=S-P$ and retains one
two-parameter pole predicate.  The same comparison point $H$ serves
the ordinary and scalar certificates.

Finite multipliers, rational and quadratic-field square roots, elliptic
ordinates, point-chart choices, and every nonzero-denominator test are
consolidated in one final rational unknown using a separable-algebra
construction.  The rational relation-combining theorem of Zhang and Sun
\cite[Theorem 1.2]{ZS} gives one explicit packaging of the square
conditions used in this last step; the full elimination is supplied by
the finite-algebra compiler of \cite[Sections 10.3--10.6]{Sun86}.
The resulting count of unknowns in our existential definition of $\Z$ in $\Q$ is
\begin{equation}\label{eq:countintro}
 \underbrace{2(4+5)}_{18\text{ point coordinates}}
 +\underbrace{1}_{\text{shared slope}}
 +\underbrace{2\cdot2}_{\text{ordinary poles}}
 +\underbrace{2}_{\text{scalar pole}}
 +\underbrace{1}_{\text{algebraic compiler}}=26.
\end{equation}

For any number field $K$ (which is a finite extension of the field $\Q$), it is known that
$\Z$ is diophantine over the ring $O_K$
of algebraic integers in $K$ (cf. \cite{ABHS,KP}). Here we pose the following general conjecture.

\begin{conjecture} Let $K$ be any number field. Then, there is a polynomial
$P(x_0,x_1,\ldots,x_n)$ with integer coefficients, such that $x_0\in K$ lies in $\Z$
if and only if $P(x_0,x_1,\ldots,x_n)=0$ for some $x_1,\ldots,x_n\in K$.
\end{conjecture}

In view of Sun's 11 unknowns theorem (cf. \cite{Sun11}), this conjecture implies that for each number field $K$,
there is a positive number $m$ such that there is no algorithm to decide for any polynomial 
$P(x_1,\ldots,x_m)\in \Z[x_1,\ldots,x_m]$ whether  $P(x_1,\ldots,x_m)=0$ for some $x_1,\ldots,x_m\in K$.

\section{Notation and arithmetic inputs}\label{sec:inputs}

The conventions in this section apply throughout the proof.  We put
$\mathbb Q^\times=\mathbb Q\setminus\{0\}$, and for any field
$L$ we write $L^\times=L\setminus\{0\}$.  The letter $p$ denotes a
rational prime.  We write $\mathbb Q_p$, $\mathbb Z_p$ and
$\mathbb F_p$ for the $p$-adic field, its valuation ring and its residue
field.  We write $\mathbb P^1(\mathbb Q)$ for the rational
projective line.  The normalized additive valuation on $\mathbb Q$ is
$v_p(p)=1$, with $v_p(0)=+\infty$.  If $L$ is a number field, then
$\mathcal O_L$ denotes its ring of integers.  For a nonzero prime ideal
$\mathcal P\mid p$ of $\mathcal O_L$, let
$\ord_{\mathcal P}(a)$ be the exponent of $\mathcal P$ in the
principal fractional ideal $(a)$ and put
$e_{\mathcal P}=\ord_{\mathcal P}(p)$.  Define
\[v_{\mathcal P}(0)=+\infty,\ \ \t{and}\ \
 v_{\mathcal P}(a)=\frac{\ord_{\mathcal P}(a)}{e_{\mathcal P}}
 \ \ \t{for}\ a\in L^\times.
\]
Thus $v_{\mathcal P}(p)=1$ and the restriction of
$v_{\mathcal P}$ to $\mathbb Q$ is $v_p$.  An element is called a
$p$-unit, or a $\mathcal P$-unit, when its corresponding valuation is
zero.  We write
$L_{\mathcal P}$ and $\mathcal O_{\mathcal P}$ for the corresponding completion
and valuation ring.  A congruence modulo $\mathcal P^k$ in this
manuscript means that the normalized valuation of the difference is at
least $k$; it is used only at unramified primes when an integral
exponent is required.

For $a\in\mathbb Q^\times$, let
\[
 \Supp(a)=\{p:v_p(a)\ne0\}.
\]
For $r=a/b\in\mathbb Q$ in lowest terms with $b>0$, we put
$\num(r)=a$.  Greatest common divisors and least common multiples of integers
are taken to be nonnegative.  We write
\[
 \Sq(\mathbb Q)=\{a^2:a\in\mathbb Q\},\qquad
 \mathbb Q^{\times2}=\Sq(\mathbb Q)\setminus\{0\},
\]
and
\[
 \mathscr D_2=\{b\in\mathbb Q:\text{the positive reduced denominator
 of $b$ is a square}\}.
\]
For a field $L$ of characteristic different from two, $[a]_L$ denotes
the class of $a\in L^\times$ in $L^\times/L^{\times2}$; the subscript
is omitted when the field is clear.  The trace and norm of a finite
extension $L/\mathbb Q$ are denoted by $\Tr_{L/\mathbb Q}$ and
$\Nm_{L/\mathbb Q}$.

The symbol $\hW$ denotes the absolute logarithmic Weil height;
in the cited sources it is denoted simply by $h$.  If
$L$ contains $\alpha$, let $\mathcal V_\infty(L)$ be the set of
archimedean places of $L$.  For every nonzero prime ideal
$\mathcal P\mid p$ of $\mathcal O_L$, put
$n_{\mathcal P}=[L_{\mathcal P}:\mathbb Q_p]$ and
$|a|_{\mathcal P}=p^{-v_{\mathcal P}(a)}$.  For
$v\in\mathcal V_\infty(L)$, put $n_v=[L_v:\mathbb R]$ and use the
ordinary real or complex modulus $|\cdot|_v$.  Thus, with the
product-formula normalization,
\[
\begin{aligned}
 \hW(\alpha)=\frac1{[L:\mathbb Q]}\Bigg(&
 \sum_{\substack{\mathcal P\subset\mathcal O_L\\\mathcal P\ne0}}
 n_{\mathcal P}\log\max\{1,|\alpha|_{\mathcal P}\}
 +\sum_{v\in\mathcal V_\infty(L)}
 n_v\log\max\{1,|\alpha|_v\}\Bigg).
\end{aligned}
\]
This quantity is independent of the chosen field $L$.  In particular,
if $a,b\in\mathbb Z$ are coprime and $b>0$, then
$ \hW(a/b)=\log\max\{|a|,b\}.$ Note that $\hW(0)=\hW(0/1)=0$.
For a rational projective point $(u:v)$ represented by coprime integers,
we use the compatible convention
\[
 \hW(u:v)=\log\max\{|u|,|v|\};
\]
thus the affine height above is the projective height of $(a:b)$.
We use the standard inequalities for sums, products and quotients of
heights.  The Néron--Tate height on the fixed elliptic curve below is
written $\widehat h_E$; thus
$\hW(x([m]P))=2m^2\widehat h_E(P)+O_E(1)$ for $m\ne0$.
All logarithms are natural logarithms.

Write $K=\mathbb Q(\rho)$, where $\rho^2=2$, and let $\sigma$ be the
nontrivial automorphism of $K/\mathbb Q$.  Thus
$(a+b\rho)^\sigma=a-b\rho$.  We put $\mathcal O_K=\mathbb Z[\rho]$.
For a rational prime $p$, we identify
$\mathcal O_K\otimes_{\mathbb Z}\mathbb Z_p$ with
$\mathbb Z_p[\rho]$ through the basis $1,\rho$; membership in
$p(\mathcal O_K\otimes\mathbb Z_p)$ therefore means coefficientwise
divisibility by $p$.  A superscript $\sigma$ on a curve, point,
coefficient or expression
means coefficientwise conjugation.  For an elliptic curve, $O$ denotes the identity, $[n]Q$ denotes
scalar multiplication by $n\in\mathbb Z$, $E[2](K)$ denotes the group
of $K$-rational points killed by $2$, and $E(K)_{\rm tors}$ denotes the
torsion subgroup.  We write $$2E(K)=\{[2]R:R\in E(K)\}\ \ \t{and}\ \ 
\mathbb ZP=\{[m]P:m\in\mathbb Z\}.$$  The symbols $x(Q)$ and $y(Q)$
denote the affine coordinates when $Q\ne O$.  A point is said to be
\emph{formal at} a nonzero prime ideal $\mathcal P$ when its image in
the corresponding completion lies in the formal neighborhood of the
identity for the integral Weierstrass equation being used.  At the good
unramified prime ideals in the contact argument this is the kernel of
reduction; at the finitely many other prime ideals we work inside an
explicitly chosen sufficiently deep open formal subgroup.
Closures in real Lie groups are denoted by an overline, and the
superscript ${}^0$ denotes the identity component.

The notation $O_{\mathcal A}(g)$ means that the implied constant may
depend only on the fixed data $\mathcal A$ indicated in the subscript;
when no subscript is written, it depends only on the fixed arithmetic
data of the construction.  Every set called ``fixed'' is chosen before
the rational input and before the quantified witnesses.

For ease of reference, the principal global symbols used later are
summarized here.  The point $P$ is the fixed generator.  The points
$Q,N,S,T,Z$ are the auxiliary and four root points; $R_{\rm sc}=S-P$
is the scalar point; $H$ is the common comparison point; and $U_i,V_i$
are its primary and translated comparison points.  Their integer labels
are written $m_Q,n,s,t,z,h,m_i$ as appropriate; these labels are used in
the proof and are never additional rational unknowns.  The functions
$\varrho(V)=1/x(V)$ and $z(V)=-x(V)/y(V)$, with their values at $O$
defined below, are local point parameters.  The symbols $\ell$ and
$\ell_{64}$ denote the formal logarithm and its degree-$64$ truncation.
The relation $\Phi_0$ is the two-structural-parameter pole predicate.
The symbols $\Delta$, $J$ and $B$ later denote, respectively, the Lucas
error, the reduced numerator of $\beta-n$, and the maximum of the four
root-index magnitudes.  Each of these symbols is formally defined again
at the point where it enters an argument.

\subsection{The cyclic carrier and the height estimate}

Put
\begin{equation}\label{eq:curve}
 d=3-\rho,\qquad D=d^2=11-6\rho,\qquad
 E:\quad y^2=x^3-Dx.
\end{equation}
The norm of $d$ is $d d^\sigma=7$, so $(d)$ is a prime ideal above
$7$. For
\[
 P=\left(\frac{57}{8}-\frac72\rho,
                     -\frac{165}{8}+\frac{417}{32}\rho\right),\]
 we have 
$[2]R_0=P$, where $R_0=(1+\rho,2\rho)$.
We write $O$ for the identity point. Burungale \cite[Proposition 2.1]{Burungale} proves that this curve has
rank one and torsion subgroup $E[2](K)$. For the particular displayed
point $P$ we use the generator calculation in
\cite[Proposition 3.1]{Glazer}, which gives
\begin{equation}\label{eq:carrier}
 2E(K)=\mathbb ZP,\qquad E(K)_{\rm tors}=E[2](K),
 \ \t{and}\  P\text{ has infinite order}.
\end{equation}
For comparison with the formal-group model in \cite{Sun86}, we set
\[
 E_0:y^2=x^3-(75-53\rho)^2x\ \ \t{and}\ \ 
 P_0=\left(\frac{1641-1160\rho}{8},
             \frac{-123720+87483\rho}{32}\right).
\]
For $\lambda_0=3-2\rho$, clearly
$(3-\rho)\lambda_0^2=75-53\rho$, and
\begin{equation}\label{eq:iso}
 (x,y)\longmapsto(\lambda_0^2x,\lambda_0^3y)
\end{equation}
is an isomorphism $E\to E_0$ carrying $P$ to $P_0$.
The element $\lambda_0$ is a unit of $\mathcal O_K$, so this integral
change of variables transports the formal-group statements used below
without changing normalized local depths outside the same fixed finite
set. For a
point $V\in\mathbb ZP$, we define
\begin{equation}\label{eq:rho}
 \varrho(V)=\begin{cases}1/x(V),&\t{if}\ V\ne O,\\0,&\t{if}\ V=O,\end{cases}
\qquad
 z(V)=\begin{cases}-x(V)/y(V),&\t{if}\ V\ne O,\\0,&\t{if}\ V=O.\end{cases}
\end{equation}
The nonidentity points of $\mathbb ZP$ have neither $x=0$ nor $y=0$:
such a point would be nonzero two-torsion, contrary to
\eqref{eq:carrier}. Values at $O$ are the regular formal continuations.

Let $\zeta=-x/y$ be an indeterminate for the formal parameter at
$O$. The invariant differential is
$(1+4D\zeta^4)^{-1/2}d\zeta$, so the formal logarithm is
\begin{equation}\label{eq:formallog}
\begin{aligned}
 \ell(\zeta)&=\sum_{j\ge0}(-1)^j\binom{2j}{j}
       \frac{D^j}{4j+1}\zeta^{4j+1},\\
 \ell_{64}(\zeta)&=\sum_{4j+1\le64}(-1)^j\binom{2j}{j}
       \frac{D^j}{4j+1}\zeta^{4j+1}.
\end{aligned}
\end{equation}
The conjugate series is obtained by applying $\sigma$ to the
coefficients. Proposition 2.3 of \cite{Burungale} proves the two facts
used below: at a good prime the valuation of the quotient of two formal
parameters equals the valuation of the quotient of their integer
indices, and for $p>65$ the quotient of the truncated logarithms agrees
with that index quotient to precision $64$ whenever the latter is
integral. The proof of that proposition uses only that $P$ generates
$2E(K)$, so it applies to the generator fixed in \eqref{eq:carrier}.
These are also the estimates used in \cite[Lemma 2.1]{Sun86}. The lower-case symbol $z(V)$ always means the
point parameter in \eqref{eq:rho}; the integer $z$ introduced later is
an index and is a different object.

The five-contact height estimate needed below follows from Burungale's
Theorem E \cite{Burungale}. His branch values are
\[
 \mathscr B=\{0,13,525,637,650\}.
\]
Consider the fractional-linear transformation
\begin{equation}\label{eq:burungale-mobius}
 \chi(s)=\frac{49s-31213}{s-13}.
\end{equation}
A direct calculation gives
\[
 \chi(13)=\infty,\quad \chi(637)=0,\quad
 \chi(650)=1,\quad \chi(0)=2401,\quad
 \chi(525)=-\frac{343}{32}.
\]
Thus the image of $\mathscr B$ is the contact set
\begin{equation}\label{eq:contactold}
 \mathscr C=\left\{\infty,0,1,2401,-\frac{343}{32}\right\}
 \subset\mathbb P^1(\mathbb Q).
\end{equation}
For $c\in\mathbb Q$ choose a primitive integral linear form
$L_c(U,V)$ vanishing at $c$; ``primitive'' means that its two integer
coefficients are coprime. Put $L_\infty(U,V)=V$. For
$a=(u:v)\in\mathbb P^1(\mathbb Q)$ and a rational prime $p$, define
its normalized contact depth at $c$ by
\begin{equation}\label{eq:contactdepth}
 d_{p,c}(a)=v_p(L_c(u,v))-\min\{v_p(u),v_p(v)\}.
\end{equation}
This does not depend on the homogeneous representative $(u:v)$. For a
fixed finite set $\Sigma$ of rational primes define
\begin{equation}\label{eq:M5S}
 M_{5,\Sigma}(a)=
 \prod_{\substack{p\notin\Sigma\\
       d_{p,c}(a)>0\text{ is odd for some }c\in\mathscr C}}p.
\end{equation}
Theorem E of \cite{Burungale}, transported by the fixed map
\eqref{eq:burungale-mobius}, gives a fixed finite exceptional set
$\Sigma_*$ and a constant $C_5>0$ such that
\begin{equation}\label{eq:fiveheight}
 \hW(a)\le C_5M_{5,\Sigma_*}(a)^{24}
 \qquad(a\in\mathbb P^1(\mathbb Q)\setminus\mathscr C).
\end{equation}
Indeed, after enlarging the exceptional set by the primes dividing the
coefficients and determinant of $\chi$, the pullbacks of primitive
contact forms differ only by local units, so odd contact is preserved;
a fixed fractional-linear transformation changes logarithmic height by
$O(1)$. Enlarging $\Sigma_*$ further by a fixed finite set only changes
the constant by a fixed factor. We absorb all such fixed changes below
into constants called $C_{\rm ht}$ or $C_{\rm sc}$, each at least one.

\subsection{Odd-index certificates and aggregate trace height}
\label{subsec:oddproduct}

Put $\varepsilon_K=1+\rho$ and $T_0=(0,0)$, $T_d=(d,0)$.  The following
finite certificate is the form of \cite[Proposition J.7]{Glazer}
needed below.

\begin{lemma}\label{lem:oddchart}
Let $Q=(x,y)\in2E(K)$ be affine.  Its unique index in
$2E(K)=\mathbb ZP$ is odd if and only if there are
$h,j,u,v\in K$ with $h\ne0$ such that
\begin{equation}\label{eq:oddchart}
 h^2=x,\qquad y=hj,\qquad
 u^2=2(x+j),\qquad uh+x+j=\varepsilon_K v^2.
\end{equation}
For each fixed $Q$, the set of such quadruples is finite.
\end{lemma}
\begin{proof}
For $A=(x_A,y_A)\in E(K)$ with $x_A\ne0$, let
$\alpha(A)$ be the squareclass $[x_A]_K$ in $K^\times/K^{\times2}$; put
$\alpha(O)=1$ and $\alpha(T_0)=[-D]_K=[-1]_K$.  The usual
line-intersection identity, together with the duplication and
$T_0$-translation formulas, shows that $\alpha$ is a homomorphism.
By \eqref{eq:carrier}, every point $A\in E(K)$ is uniquely
$mR_0+aT_0+bT_d$, with $m\in\mathbb Z$ and $a,b\in\{0,1\}$.
Indeed, $[2]A=kP=[2k]R_0$ for a unique $k\in\mathbb Z$, so
$A-kR_0\in E[2](K)$; uniqueness follows from the infinite order of
$R_0$.  Hence
\[
 \alpha(mR_0+aT_0+bT_d)=[\varepsilon_K]_K^m[-1]_K^a[d]_K^b.
\]
The classes $[\varepsilon_K]_K,[-1]_K,[d]_K$ are independent.  Indeed,
for $a,b,c\in\{0,1\}$ the norm of
$\varepsilon_K^c(-1)^a d^b$ is $(-1)^c7^b$; if the element is a square in
$K$, its norm is a nonnegative rational square, so $b=c=0$, and the
two real embeddings then force $a=0$.
Consequently
\begin{equation}\label{eq:alphacoset}
 \alpha(A)=[\varepsilon_K]_K\quad\Longleftrightarrow\quad
 A\in R_0+2E(K).
\end{equation}

Assume \eqref{eq:oddchart}.  The curve equation and $h\ne0$ give
$j^2=x^2-D$.  Put
\[
 w=u^2=2(x+j),\qquad X_*=uh+\frac w2=uh+x+j,
 \qquad Y_*=X_*u.
\]
Then
\[
 w\left(x-\frac w4\right)=D,
 \qquad X_*^2-wX_*=D,
 \qquad Y_*^2=X_*^3-DX_*.
\]
Since $D\ne0$, the first two identities give $u\ne0$ and
$X_*\ne0$; hence $Y_*=X_*u\ne0$.  Thus $(X_*,Y_*)\in E(K)$ is a
nontorsion affine point, and direct substitution in the legitimate tangent-duplication
formulas shows that $[2](X_*,Y_*)=Q$.  The last equation in
\eqref{eq:oddchart} says $X_*=\varepsilon_K v^2$, so
\eqref{eq:alphacoset} gives $(X_*,Y_*)\in R_0+2E(K)$.  Hence
$Q\in P+4E(K)$, which is precisely the set of odd multiples of $P$.

Conversely, suppose $Q=[m]P$ with $m$ odd.  Take
$(X_*,Y_*)=[m]R_0$, so $[2](X_*,Y_*)=Q$ and
$\alpha((X_*,Y_*))=[\varepsilon_K]_K$.  Thus $X_*=\varepsilon_K v^2$ for
some $v\in K$.  Set
\[
 u=Y_*/X_*,\qquad h=(X_*^2+D)/(2Y_*),
 \qquad x=h^2,\qquad j=u^2/2-x.
\]
The duplication formulas give $Q=(x,hj)$, and direct substitution
gives both remaining equations in \eqref{eq:oddchart}.  Since the odd
multiple $Q$ is nontorsion, $x(Q)\ne0$, and therefore the produced
$h$ is nonzero.  The choices of $h,u,v$ are roots of three fixed
quadratic equations over the coordinates of $Q$, and $j=y/h$ is then
unique. Hence the fiber is finite.
\end{proof}

\begin{lemma}\label{lem:aggregatetrace}
For nonzero integers $m_1,\ldots,m_4$, put
\begin{equation}\label{eq:aggregate-trace}
 T_{\rm ag}(m_1,\ldots,m_4)
 =\frac12\Tr_{K/\Q}\left(\sum_{i=1}^4x([m_i]P)\right).
\end{equation}
Then $T_{\rm ag}>0$, and there are fixed constants $M_1>0$ and
$M_2\ge0$, depending only on the fixed curve and $P$, such that
\begin{equation}\label{eq:aggregate-height}
 \max_i|m_i|^2\le M_1\hW(T_{\rm ag})+M_2.
\end{equation}
No parity assumption on the $m_i$ is needed.
\end{lemma}
\begin{proof}
Every nonzero point $[m]P$ is a nonzero double, since $P=[2]R_0$.
For the curve $y^2=x^3-Dx$, the duplication formula gives
\[
 x([2]A)=\left(\frac{x(A)^2+D}{2y(A)}\right)^2
\]
whenever $A$ is not two-torsion.  Hence $x([m]P)$ is a nonzero square
in $K$ and is positive at both real embeddings.  This proves
$T_{\rm ag}>0$.

The sum version of the aggregate trace-height theorem in
\cite[Section A38]{Glazer}, applied
with $r=4$ to this same curve and generator, gives fixed constants
$M_+>0$, $M_{+,0}\ge0$ such that
\[
 \sum_{i=1}^4\widehat h_E([m_i]P)
 \le M_+\hW(T_{\rm ag})+M_{+,0}.
\]
Since $P$ is nontorsion,
$\widehat h_E([m]P)=m^2\widehat h_E(P)$ and
$\widehat h_E(P)>0$.  Dividing by the fixed positive canonical height
and enlarging the additive constant gives \eqref{eq:aggregate-height}.
\end{proof}

\begin{corollary}\label{cor:twoscales}
For nonzero integers $m_1,m_2,m_3,m_4$, set
\[
 T_{\rm ag}=T_{\rm ag}(m_1,m_2,m_3,m_4),\qquad
 A_0=1\ \ \t{and}\ \ A_1=T_{\rm ag}^2.
\]
If $\eta\in\Q^\times$ and
\[
 \hW(A_j\eta)\le H_0\qquad(0\le j\le1),
\]
then
\[
 \max_i|m_i|^2\le M_1H_0+M_2.
\]
\end{corollary}
\begin{proof}
The case $j=0$ gives $\hW(\eta)\le H_0$.  The standard height
inequality gives
\[
 2\hW(T_{\rm ag})=\hW(T_{\rm ag}^2)
 \le\hW(T_{\rm ag}^2\eta)+\hW(\eta)\le2H_0.
\]
Now apply Lemma~\ref{lem:aggregatetrace}.
\end{proof}

\subsection{The two-structural-parameter pole predicate}

We use the precise construction of \cite[Proposition 2.3]{Sun86}.
Fix the finite nonzero multiplier set $\mathcal M\subset\mathbb Q^\times$
there. Let
\[
 S_0=\{2,3\}\cup\bigcup_{\mu\in\mathcal M}\Supp(\mu),
\]
where $\Supp(\mu)=\{p:v_p(\mu)\ne0\}$. For $b\in\mathbb Q^\times$,
the relation $\Phi_0(b)$ is the finite disjunction, over $\mu\in\mathcal M$,
of the following conditions, with only $e,w\in\mathbb Q$ as free
structural parameters:
\begin{gather}
 b'=\mu b,\qquad h'=b'(b'-b^2),\label{eq:poledata}\\
 b'e h'w\ne0,\label{eq:poleguard}\\
 b'(b'-e)(b'+3e)\in\Sq(\mathbb Q),\label{eq:polesq1}\\
 b'(b'w^2-eh')(b'w^2+3eh')\in\Sq(\mathbb Q).
 \label{eq:polesq2}
\end{gather}
The extension $\Phi_0(0)$ is true.  The multiplier $\mu$ is selected
from a fixed finite list and is therefore a finite branch label.  The
two square roots are finite algebraic choices.  Neither contributes a
continuous structural parameter.  Proposition 2.3 of \cite{Sun86}
provides a fixed constant $\varepsilon>0$ and fixed integers $b_p>0$
for $p\in S_0$ such that
\begin{align}
 &\Phi_0(b),\quad b\ne0\ \Longrightarrow\
     v_p(b)\notin\{-1,-3,-5,\ldots\}
       \quad(p\notin S_0),\label{eq:polesound}\\
 &
 b\in\mathscr D_2,\ 0<b<\varepsilon,
 \ v_p(b)\ge b_p\quad(p\in S_0)
 \ \Longrightarrow\ \Phi_0(b).
 \label{eq:polecomplete}
\end{align}
In the first statement a negative odd valuation is forbidden; in the second,
there is no additional local-square requirement on $b$. Both statements
are exactly the cited proposition.

Burungale \cite[Theorem B and Proposition 3.14]{Burungale} gives an
independent unconditional construction of a square-denominator
predicate using congruent-number curves, the Tian--Yuan--Zhang
analytic-rank criterion, Gross--Zagier and Kolyvagin. This shows that
the square-denominator mechanism itself need not use the pointwise
$2$-converse appearing in some earlier approaches. We do not replace
$\Phi_0$ by Burungale's alternative here, because the predicate
\eqref{eq:poledata}--\eqref{eq:polesq2} is already integrated into the
local preparation and has exactly the two structural parameters used in
the counts in Theorems~\ref{thm:main} and \ref{thm:abc16}.

\subsection{Simultaneous orders and formal logarithms}

Let $i^2=-1$.  Thus $\mathbb Z[i]$ is the ring of Gaussian integers
and $\mathbb Q(i)$ its fraction field.  For a rational prime $p$ and a
Gaussian prime $\pi=a+bi\in\mathbb Z[i]$ of norm $p$, put
\begin{align*}
 d_+(p)&=\Nm_{\mathbb Q(i)/\mathbb Q}(\pi+1)=(a+1)^2+b^2,\\
 d_-(p)&=\Nm_{\mathbb Q(i)/\mathbb Q}(\pi-1)=(a-1)^2+b^2.
\end{align*}
For the primes supplied below, these are the cardinalities of the two
reductions of $E$ at the primes of $K$ above $p$, with the orientation
fixed by the Gaussian congruences.  They annihilate the corresponding
reduced points; exact orders of the reduction of $P$ are not required.

The special instance of \cite[Theorem 4.1]{Sun86} that we require is
the following.  For any fixed distinct integral shifts
$c_1,\ldots,c_k$ ($k\ge2$), any positive integer $q_{\rm div}$, any odd
squarefree positive $M_0$ not divisible by $7$, and any positive lower
bounds, there exist positive integers $u,w$ with $q_{\rm div}\mid w$ and
$u/w^2$ arbitrarily large such that every
\[
 p_j=u+c_jw^2
\]
is a distinct prime.  Each $p_j$ has a Gaussian prime generator
$\pi_j$ satisfying the following congruences in $\mathbb Z[i]$:
\[
 \pi_j\equiv9\pmod {16},\qquad
 \pi_j\equiv1+3i\pmod 7,\qquad
 \pi_j\equiv i\pmod {M_0}.
\]
In particular,
\begin{gather}
 \gcd(d_+(p_i),d_-(p_j))=4\quad(1\le i,j\le k),
 \label{eq:crossorders}\\
 v_2(d_+(p_j))=2,\quad v_2(d_-(p_j))=6,
 \quad \gcd(d_+(p_j)d_-(p_j),M_0)=1.
 \label{eq:ordproperties}
\end{gather}
All these cardinalities are prime to $7$, by the fixed Gaussian congruence
modulo $7$ in that theorem. Its formulation also permits a prescribed
square-denominator dilation. Repeated shift positions are removed
before invoking the theorem. The requirement that the shifts be
fixed before $u$ and $w$ are chosen is essential.

Let $c_*$ be the fixed positive integer in the two-anchor construction
of \cite[Sections 2.5 and 3]{Sun86}, transported through
\eqref{eq:iso}.  It is chosen so that, for every $p\in S_0$ and every
$\mathcal P\mid p$, the point $[c_*]P$ lies in the formal subgroup
with $v_{\mathcal P}(z([c_*]P))\ge2$.  Enlarging $c_*$ by a fixed
multiple preserves this property, so we also require $64\mid c_*$.  We
use the truncations
$\ell_{64}$ and $\ell_{64}^{\sigma}$ defined in
\eqref{eq:formallog}. Sun's local comparison lemma, transported across
\eqref{eq:iso}, applies to these logarithms at primes
$p>65$ outside a fixed exceptional set. The precise four-point
application is proved in Section~\ref{sec:ordinary}.

\subsection{A rational relation-combining lemma}

We need the following useful result of Zhang and Sun \cite{ZS}.

\begin{lemma}[Zhang--Sun \cite{ZS}, Theorem 1.2]\label{lem:ZS}
For every integer $k\ge1$ there is an explicitly constructible
integer polynomial $J_k(A_1,\ldots,A_k,W)$ with the property that, for
all nonzero $A_1,\ldots,A_k\in\mathbb Q$,
\begin{equation}\label{eq:ZS}
 A_1,\ldots,A_k\in\Sq(\mathbb Q)
 \iff
 \exists w\in\mathbb Q\,[J_k(A_1,\ldots,A_k,w)=0].
\end{equation}
\end{lemma}
The full explicit product defining $J_k$ is given in the cited paper;
we use this lemma without reproving it. Square conditions with zero
input will be handled in separate finite branches. In the final
compiler, this lemma is used together with the stronger separable
primitive-element mechanism of
\cite[Sections 10.3--10.6]{Sun86}, which also absorbs other finite
algebraic choices.

\section{The ordinary sixteen-test relation}\label{sec:ordinary}

Fix a prospective rational input $X$ and set
\begin{equation}\label{eq:beta}
 \beta=2X^2+3.
\end{equation}
Take nine cyclic points
\begin{equation}\label{eq:ninepoints}
 Q,N,T,Z,H,U_1,U_2,U_3,U_4\in\mathbb ZP.
\end{equation}
A sufficiently large fixed positive integer $M$, divisible by $6$,
will be chosen in Section~\ref{sec:heights}.  Put
\begin{equation}\label{eq:Sdef}
 R_{\rm sc}=[M]Q,\qquad S=P+R_{\rm sc},
 \qquad P_1=N,\quad P_2=S,\quad P_3=T,\quad P_4=Z.
\end{equation}
Write
\[
 (P_1,P_2,P_3,P_4)=([n]P,[s]P,[t]P,[z]P)
\]
for their unique integer indices.  These indices are used only to
prove the equivalence; they are never quantified rationally.  If
$Q=[m_Q]P$, then
\begin{equation}\label{eq:scalarindex}
 s-1=Mm_Q,
 \qquad R_{\rm sc}=[s-1]P.
\end{equation}
We impose
$Q,R_{\rm sc},S,N,T,Z\ne O$ and impose the certificate of
Lemma~\ref{lem:oddchart} on $N,T,Z$.  Since $6\mid M$, the index $s$
is automatically odd.  Thus all four integers $n,s,t,z$ are odd,
$|s-1|\ge M$, and $z\ne0$.

Put $V_0=H-[c_*]P$ and
$V_i=U_i-[c_*]P_i$ for $1\le i\le4$. We impose
$H\ne O$, $V_0\ne O$, and require the following two denominators in
$K$ to be nonzero:
\begin{equation}\label{eq:logden}
 \ell_{64}(z(H)),\qquad
 \ell_{64}^{\sigma}(z(V_0^{\sigma})).
\end{equation}
Define the computed values
\begin{equation}\label{eq:logs}
 L_i=\frac{\ell_{64}(z(U_i))}{\ell_{64}(z(H))},\qquad
 L'_i=\frac{\ell_{64}^{\sigma}(z(V_i^{\sigma}))}
               {\ell_{64}^{\sigma}(z(V_0^{\sigma}))}
 \quad (1\le i\le4).
\end{equation}
Here all expressions lie in $K$; they are not new quantified variables.
Set
\begin{equation}\label{eq:lucaspoly}
 \mathcal D(X_1,X_2,X_3,X_4)
 =X_3^2-X_2(X_2-2)X_3X_4+(X_2-2)^2X_4^2-1.
\end{equation}
In particular, $\Delta=\mathcal D(n,s,t,z)$.

We enumerate sixteen $K$-valued tests $D_i$, and their required
local precisions $k_i$, in the displayed order as follows:
\begin{equation}\label{eq:sixteentests}
\begin{array}{c|c|c}
 \text{test expressions} & \text{number} & k_i\\ \hline
 \varrho(H),\ \varrho(U_1),\ldots,\varrho(U_4)&5&1\\
 \varrho(V_0^{\sigma}),\ \varrho(V_1^{\sigma}),\ldots,
      \varrho(V_4^{\sigma})&5&1\\
 L_i-L'_i\ (1\le i\le4)&4&64\\
 \mathcal D(L_1,L_2,L_3,L_4)&1&64\\
 L_1-\beta&1&16\\ \hline
 \text{total}&16&
\end{array}
\end{equation}
All denominators used in the rational group law, $\varrho$, $z$,
$\ell_{64}$, and the normalizing expressions below are required to
be nonzero on the chosen finite charts.

\begin{lemma}\label{lem:tw_anchor}
There is a fixed finite prime set $S_{\rm ord}$ containing the bad and
ramified primes, the primes dividing $c_*$, and all primes at most $65$,
with the following property. Suppose $p\notin S_{\rm ord}$ and there is
$\mathcal P\mid p$ for which
\begin{equation}\label{eq:prec_conditions}
 v_{\mathcal P}(D_i)\ge k_i\qquad(1\le i\le16).
\end{equation}
Then $v_p(\Delta)\ge64$ and
\begin{equation}\label{eq:betaprec}
 v_p(\beta-n)\ge16.
\end{equation}
In particular, if $\Delta\ne0$, then $p^{64}\mid\Delta$; if
$\beta\ne n$ and $J=|\num(\beta-n)|$ in lowest terms,
then $p^{16}\mid J$.
\end{lemma}
\begin{proof}
The first ten inequalities make $H,U_i$ formal at
$\mathcal P$ and $V_0,V_i$ formal at the conjugate prime ideal $\mathcal P^\sigma$. If a
numerator is zero, the corresponding truncated-logarithm numerator is
zero and the argument below is unchanged. The formal-logarithm estimate \cite[Proposition 2.3]{Burungale},
used in the same two-anchor arrangement as the proof of
\cite[Lemma 3.2]{Sun86}, shows that
\[
 L_i\equiv\frac{m_i}{h}\ \ \t{and}\ \ 
 L'_i\equiv\frac{m_i-c_*n_i}{h-c_*}\pmod{\mathcal P^{64}}
\]
whenever the corresponding rational ratios are integral; here
$H=[h]P$, $U_i=[m_i]P$, and
$(n_1,n_2,n_3,n_4)=(n,s,t,z)$.
Since $p\nmid c_*$, at least one of $h,h-c_*$ is a $p$-unit.
For that denominator the corresponding integer index ratio is
$p$-integral, so \cite[Proposition 2.3]{Burungale} makes the corresponding
truncated ratio integral.  The inequality
$v_{\mathcal P}(L_i-L'_i)\ge64$ then makes the other truncated
ratio integral as well, and the valuation identity in that proposition
makes the other index ratio integral. Multiplying by
$h$ and $h-c_*$ and subtracting gives
\[
 c_*L_i\equiv hL_i-(h-c_*)L'_i
 \equiv c_*n_i\pmod{\mathcal P^{64}}.
\]
Thus $L_i\equiv n_i\pmod{\mathcal P^{64}}$ for every $i$.
The polynomial $\mathcal D$ has integer coefficients, so the
next-to-last test yields $\Delta\equiv0\pmod{\mathcal P^{64}}$.
The last test yields $\beta-n\equiv0\pmod{\mathcal P^{16}}$.
Both quantities are rational, and the first is an integer;
normalization $v_{\mathcal P}(p)=1$ gives the asserted rational
valuations. The statement about the reduced numerator follows because
positive $p$-adic valuation forces the denominator to be a $p$-unit.
\end{proof}

\section{Two ordinary contact certificates}\label{sec:contacts}

Move the original infinity contact to a finite point using
\[
 \varphi(v)=2+\frac1v.
\]
The inverse image of \eqref{eq:contactold} is exactly
\begin{equation}\label{eq:thetacontacts}
 \mathscr T=\left\{0,-\frac12,-1,\frac1{2399},-\frac{32}{407}\right\}.
\end{equation}
Fix $\xi=3$, which is not among these contacts, and set
\begin{equation}\label{eq:Ffive}
 F_5(V)=\prod_{\theta\in\mathscr T}(V-\theta).
\end{equation}
For each $\theta\in\mathscr T$ choose a primitive integral linear
form $L_\theta(U,V)$ vanishing at $\theta$.  Fix once and for all a
finite set $S_{\rm ht}$ containing the transported exceptional set,
$S_0$, $S_{\rm ord}$, the coefficient primes of these forms, and the
fixed contact-difference primes.  For $\tau=(u:w)\notin\mathscr T$ put
\begin{equation}\label{eq:Mfinite}
 \mathfrak M(\tau)=
 \prod_{\substack{p\notin S_{\rm ht}\\
 v_p(L_\theta(u,w))-\min\{v_p(u),v_p(w)\}>0
 \text{ is odd for some }\theta\in\mathscr T}}p.
\end{equation}
Lemma~\ref{lem:contacttransport} below shows that the fixed projective change turns \eqref{eq:fiveheight} into
\begin{equation}\label{eq:finiteheight}
 \hW(\tau)\le C_{\rm ht}\,\mathfrak M(\tau)^{24}\ \ \t{with}\ \ 
 \tau\in\mathbb Q\setminus\mathscr T.
\end{equation}
The extra exclusion $\tau\ne\xi$ used below belongs to the detector
normalization, not to the height theorem.

\begin{lemma}
\label{lem:contacttransport}
After enlarging the exceptional set by a fixed finite set, the map
$\varphi(\tau)=2+1/\tau$ carries the normalized contact depths at
$\mathscr T$ to the normalized contact depths at $\mathscr C$.
Moreover $\hW(\tau)\le \hW(\varphi(\tau))+O(1)$.  Consequently
\eqref{eq:finiteheight} follows from \eqref{eq:fiveheight} with one
fixed enlarged constant $C_{\rm ht}$.
\end{lemma}
\begin{proof}
In homogeneous coordinates the map is
$(U:V)\mapsto(2U+V:U)$.  The pullback of the linear form defining the
old contact at infinity is $U$, and the pullback of a primitive linear
form defining any finite old contact is a fixed nonzero rational
multiple of a primitive form $U-\theta V$ with
$\theta\in\mathscr T$.  Outside the finitely many primes dividing
the two change-of-coordinate matrices, the scalar multiples and the primitive-form normalizations are all local units.
Thus the normalized depths, with the common minimum of the homogeneous
coordinate valuations subtracted, agree.  The standard height
inequality for a fixed fractional linear transformation
gives the final assertion.
\end{proof}

Put
\begin{equation}\label{eq:current-trace}
 T_{\rm ag}=\frac12\Tr_{K/\Q}
 \bigl(x(N)+x(S)+x(T)+x(Z)\bigr),
 \qquad A_0=1,\quad A_1=T_{\rm ag}^2.
\end{equation}
By Lemma~\ref{lem:aggregatetrace}, $T_{\rm ag}>0$, so both scales are
nonzero rational squares.  Introduce one free rational slope $\eta$ and set
$
 t_j=A_j\eta\qquad(0\le j\le1).$
Require $t_j\notin\mathscr T\cup\{\xi\}$.

For each $j$, define the rational detector and the normalized
$K$-valued test expressions
\begin{equation}\label{eq:Rj}
 R_j=A_j^{-5}\frac{F_5(t_j)}{(t_j-\xi)^6}\ \ \t{and}\ \ 
 \widetilde D_{j,i}=\frac{D_i}{(t_j-\xi)^{k_i+1}}.
\end{equation}
The $A_j^{-5}$ factor is a rational square. Define
\begin{align}
 \mathcal W_j&=
 \prod_{i=1}^{16}\left(1+R_j
       \left(\frac{\widetilde D_{j,i}}{R_j^{k_i}}\right)^2\right)
 \ \in K,\label{eq:Wj}\\
 \Theta_j&=\Tr_{K/\mathbb Q}(\mathcal W_j^{-1}),
 \qquad B_j=\frac{R_j}{\Theta_j^2}\in\mathbb Q.
 \label{eq:Bj}
\end{align}
Require $R_j\Theta_j\ne0$, $\mathcal W_j\ne0$, and
\begin{equation}\label{eq:ordinary_poles}
 \Phi_0(B_j)\qquad(0\le j\le1).
\end{equation}
Each copy of $\Phi_0$ requires two structural rational parameters.

\begin{lemma}\label{lem:contact}
Fix $0\le j\le1$ and let $p$ be a prime factor of
$\mathfrak M(t_j)$. Under the conditions
\eqref{eq:Wj}--\eqref{eq:ordinary_poles}, there is a prime
$\mathcal P\mid p$ such that
\[
 v_{\mathcal P}(D_i)\ge k_i\qquad(1\le i\le16).
\]
Consequently, for each $0\le j\le1$, we have
\begin{equation}\label{eq:contactdivides}
 \mathfrak M(t_j)^{64}\mid\Delta
 \quad(\Delta\ne0)\ \ \t{and}\ \ 
 \mathfrak M(t_j)^{16}\mid J
 \quad(\beta\ne n),
\end{equation}
where $J=|\num(\beta-n)|$.
\end{lemma}
\begin{proof}
At a positive finite contact, $v_p(t_j)\ge0$ and $t_j-\xi$ is a unit,
after the fixed contact-difference primes have been excluded. Hence
$v_p(R_j)$ is odd, since $A_j^{-5}$ is a square. If this odd valuation
were negative, the pole predicate itself would exclude the assignment.
The exact norm--trace implication of \cite[Lemma 3.1]{Sun86}, combined
with \eqref{eq:polesound}, therefore forces a prime ideal
$\mathcal P\mid p$ at which all normalized expressions have the stated
precision. Since $t_j-\xi$ is a local unit, the same is true of the
$D_i$. Lemma~\ref{lem:tw_anchor} gives $p^{64}\mid\Delta$ when
$\Delta\ne0$, and $p^{16}\mid J$ when $\beta\ne n$.  Multiplying
over the distinct prime factors of $\mathfrak M(t_j)$ gives
\eqref{eq:contactdivides} for this detector. Repeating the argument
for $j=0,1$ proves the two individual divisibilities; no
coprimality between different detectors is required.
\end{proof}

\section{The scalar test and its valuation theorem}\label{sec:scalar}

For $\beta$ as in \eqref{eq:beta} put
\begin{equation}\label{eq:r0}
 r_0(\beta)=\beta(\beta-1)(\beta-2401)
                     (32\beta+343)(2\beta-1).
\end{equation}
This never vanishes for $X\in\mathbb Q$.  Indeed $\beta\ge3$, and the
only possible remaining zero, $\beta=2401$, would give
$X^2=1199$.  A rational square which is an integer is an integer
square, while $34^2<1199<35^2$.

The six scalar tests are the two rational components of
\begin{equation}\label{eq:scalarexpr}
 \varrho(H),\qquad \varrho(R_{\rm sc}),\qquad
 \frac{z(R_{\rm sc})}{z(H)}
\end{equation}
in the basis $1,\rho$ of $K$. Denote them by $G_1,\ldots,G_6$, in
this order, and assign orders
\begin{equation}\label{eq:scalarorders}
 (k_1^*,\ldots,k_6^*)=(1,1,1,1,48,48).
\end{equation}
Write $r_0=r_0(\beta)$. Apply the quadratic norm operation
$N_{r_0}(a,b)=a^2+r_0b^2$ to the seven leaves
\begin{equation}\label{eq:scalarleaves}
 1,\quad\frac{G_1}{r_0},\ \frac{G_2}{r_0},\
 \frac{G_3}{r_0},\ \frac{G_4}{r_0},\
 \frac{G_5}{r_0^{48}},\ \frac{G_6}{r_0^{48}}
\end{equation}
using the binary tree
\begin{align}
 a_{\rm sc}&=1+\frac{G_1^2}{r_0},&
 b_{\rm sc}&=\frac{G_2^2}{r_0^2}+\frac{G_3^2}{r_0},\nonumber\\
 c_{\rm sc}&=\frac{G_4}{r_0},&
 d_{\rm sc}&=\frac{G_5^2}{r_0^{96}}+\frac{G_6^2}{r_0^{95}},\nonumber\\
 \mathcal W_{\rm sc}
 &=\bigl(a_{\rm sc}^2+r_0b_{\rm sc}^2\bigr)^2
   +r_0\bigl(c_{\rm sc}^2+r_0d_{\rm sc}^2\bigr)^2.
 \label{eq:scalarroot}
\end{align}
The constant leaf $1$ is on the all-left path. Require
\begin{equation}\label{eq:Wscguard}
 \mathcal W_{\rm sc}\ne0.
\end{equation}

Let $\Gamma,D_0$ be fixed positive integers, chosen in
\eqref{eq:Gamma}--\eqref{eq:Dfixed} before the input is quantified.
Take a finite sign branch $\delta\in\{\pm1\}$ and set
\begin{equation}\label{eq:Bscalar}
 r=\delta r_0,\qquad
 B_{\rm sc}=\frac{\Gamma^2r\mathcal W_{\rm sc}^2}
 {(1+D_0r^2)^2}.
\end{equation}
Since $r\ne0$ and \eqref{eq:Wscguard} holds, $B_{\rm sc}\ne0$.
The scalar certificate is
\begin{equation}\label{eq:scalarpole}
 \Phi_0(B_{\rm sc}).
\end{equation}
It uses two continuous structural parameters.

\begin{lemma}[Scalar odd-contact soundness]\label{lem:scalarsound}
There is a fixed finite prime set $S_{\rm sc}$, containing the
exceptional set in \eqref{eq:fiveheight}, such that every solution of
the scalar certificate satisfies
\begin{equation}\label{eq:scalarDiv}
 M_{5,S_{\rm sc}}(\beta)^{48}\mid(s-1).
\end{equation}
Consequently, there is a fixed $C_{\rm sc}>0$ such that
\begin{equation}\label{eq:scalarheight}
 \hW(\beta)\le C_{\rm sc}\sqrt{|s-1|}.
\end{equation}
\end{lemma}
\begin{proof}
Take $S_{\rm sc}$ to contain $S_0$, all primes at most $76176$,
$\Supp(\Gamma D_0)$, the bad and ramified primes of $E$, and the
finitely many coefficient and contact-difference primes used below.  Fix $p\notin S_{\rm sc}$.  Suppose that $p$ is a positive odd contact of $\beta$.
If $v_p(X)<0$, then $p\ne2$ and
$v_p(\beta)=2v_p(X)$ is negative and even, so no odd infinity contact
occurs. At a finite odd contact, $\beta$ is $p$-integral and exactly
one corresponding factor of \eqref{eq:r0} has positive odd valuation.
Excluding the fixed coefficient and contact-difference primes, put
\[
 R=v_p(r_0)>0,
 \qquad R\text{ odd}.
\]
Since $p\nmid\Gamma$ and $1+D_0r^2$ is a $p$-unit,
\begin{equation}\label{eq:rootval}
 v_p(B_{\rm sc})=R+2v_p(\mathcal W_{\rm sc}).
\end{equation}
This valuation is odd. By \eqref{eq:polesound} and the nonzero guard it
is positive, so
\begin{equation}\label{eq:rootlower}
 v_p(\mathcal W_{\rm sc})\ge(1-R)/2>-R/2.
\end{equation}
At every node $u^2+r_0v^2$, the two nonzero summands have opposite
valuation parity, so no cancellation occurs. If a node has valuation
larger than $-cR$, with $c<1$, then its children satisfy
\[
 v_p(u)>-cR/2\ \ \t{and}\ \ 
  v_p(v)>-(c+1)R/2.
\]
Starting with $c=1/2$ at the root and applying this rule down
the displayed tree gives the following coefficients of $R$ in the
strict lower bounds for the seven leaves:
\[
\begin{array}{c|ccccccc}
\text{leaf}&1&G_1/r_0&G_2/r_0&G_3/r_0&G_4/r_0&G_5/r_0^{48}&G_6/r_0^{48}\\ \hline
\text{coefficient}&1/16&9/16&5/16&13/16&3/8&7/16&15/16.
\end{array}
\]
Zero summands are interpreted as having valuation $+\infty$ and cause
no cancellation.  Every displayed coefficient is strictly smaller
than one.  Hence every normalized nonconstant leaf has valuation
strictly greater than $-R$, that is,
\[
 v_p(G_i/r_0^{k_i^*})>-R\qquad(1\le i\le6).
\]
As the valuations are integral,
$v_p(G_i)\ge k_i^*$. The first four inequalities put $H$ and
$R_{\rm sc}$ in the formal groups at every prime ideal $\mathcal P\mid p$; the last
two give
\[
 v_{\mathcal P}\left(\frac{z(R_{\rm sc})}{z(H)}\right)\ge48.
\]
Writing $H=[h]P$ and using \eqref{eq:scalarindex}, the group
identity
\[
 [h]R_{\rm sc}=[s-1]H
\]
gives
\[
 h\,\ell(z(R_{\rm sc}))=(s-1)\ell(z(H))
\]
in the formal group.  Since $\ell(z(V))/z(V)$ is a local unit there,
the formal-logarithm isometry gives
\[
 v_{\mathcal P}\left(\frac{z(R_{\rm sc})}{z(H)}\right)
 =v_p((s-1)/h)=v_p(s-1)-v_p(h).
\]
The nonzero integer $h$ has nonnegative $p$-valuation, so
$v_p(s-1)\ge48$. Multiplication over the odd contact primes outside $S_{\rm sc}$ proves
\eqref{eq:scalarDiv}.  The enlarged-set form of
\eqref{eq:fiveheight} then gives \eqref{eq:scalarheight}.
\end{proof}

\section{Simultaneous completeness with one comparison point}\label{sec:complete}

The quantifier order is important.  The four root points and the two
scales are fixed before the simultaneous-prime theorem is invoked.
Only after the contact primes are chosen is the common comparison index
selected.

\subsection{Real bounds independent of congruence moduli}

Let $E^\sigma$ be the conjugate real curve and put
\[
 G_{\R}=\overline{\{[m](P,P^\sigma):m\in\mathbb Z\}}
       \subseteq E(\mathbb R)\times E^\sigma(\mathbb R).
\]
This is a compact abelian Lie group; $G_{\R}^0$ denotes its identity
component.  In this subsection a point $V\in E(K)$ is identified with
$(V,V^\sigma)$ in the product of the two real curves.  Since
multiplication by a positive integer is surjective
on the torus $G_{\R}^0$, the closure of
$\{[h_0+Lk](P,P^\sigma):k\ge0\}$ contains a full coset of
$G_{\R}^0$ for every $L>0$ and $h_0\in\mathbb Z$.  The positive
tail is dense in the same closure as the two-sided progression.  This
is the compact-group argument used in \cite[Lemma 8.1(2)]{Sun86}.

\begin{lemma}\label{lem:real}
There are fixed open neighborhoods $\mathcal O_{\rm sc}$ of the
identity in $G_{\R}^0$, a fixed open set $\mathcal O_H$ meeting every
component of $G_{\R}$, and a constant $A_{\rm sc}\ge1$, independent of the
rational input and of later congruence moduli, with the following
properties. If $R_{\rm sc}=[s-1]P\in\mathcal O_{\rm sc}$ and
$H=[h]P\in\mathcal O_H$, all six rational coefficients of
\eqref{eq:scalarexpr} have absolute value at most $A_{\rm sc}$. Moreover, for
fixed root indices $(n,s,t,z)$, every nonconstant arithmetic
progression contains arbitrarily large $h$ for which $H\in\mathcal O_H$,
all sixteen expressions \eqref{eq:sixteentests} are defined, and their
two real embeddings are bounded by a constant depending only on the
fixed roots and arithmetic data, not on the progression modulus.
\end{lemma}
\begin{proof}
On a real identity neighborhood, $\varrho(R_{\rm sc})$ and
$z(R_{\rm sc})$ tend to zero. Choose $\mathcal O_{\rm sc}$ small
and choose $\mathcal O_H$ away from the finitely many zeros and poles
of $\varrho(H)$ and $z(H)$, and away from $H=O$. Then
$\varrho(H)$, $\varrho(R_{\rm sc})$, and
$z(R_{\rm sc})/z(H)$ are uniformly bounded at both real embeddings.
Passing to rational coefficients in the basis $1,\rho$ is a fixed
linear operation.

For fixed nonzero target indices, put $U_i=[hn_i]P$ as in the
completeness construction.  The required conclusion for the ordinary
expressions is the fixed-target compactness statement proved in
\cite[Lemma 8.1(2)]{Sun86}, transported through the integral
isomorphism \eqref{eq:iso}.  Concretely, the complete finite chart
atlas removes the finitely many denominator-zero loci, leaving a
nonempty open subset of each relevant coset on which all rational
functions are bounded.  The preceding density statement makes every
nonconstant arithmetic progression meet that subset arbitrarily far
out.
\end{proof}

In the completeness argument the input $X$ is integral, so
$r_0(\beta)\in\mathbb Z\setminus\{0\}$.  Put
$R_{\rm abs}=|r_0(\beta)|\ge1$.  If $|G_i|\le A_{\rm sc}$, then
\[
 |a_{\rm sc}|\le1+A_{\rm sc}^2\ \ \t{and}\ \ 
 |b_{\rm sc}|\le 2A_{\rm sc}^2/R_{\rm abs},
\]
so
$|a_{\rm sc}^2+r_0b_{\rm sc}^2|\le(1+A_{\rm sc}^2)^2+4A_{\rm sc}^4$.  Likewise
\[
 |c_{\rm sc}|\le A_{\rm sc}/R_{\rm abs}\ \ \t{and}\ \ 
 |d_{\rm sc}|\le2A_{\rm sc}^2/R_{\rm abs}^{95},
\]
whence
$|c_{\rm sc}^2+r_0d_{\rm sc}^2|\le(A_{\rm sc}^2+4A_{\rm sc}^4)/R_{\rm abs}^2$.  Substitution in the last norm node
therefore gives the fixed bound
\begin{equation}\label{eq:CW}
 |\mathcal W_{\rm sc}|\le C_W:=\bigl((1+A_{\rm sc}^2)^2+4A_{\rm sc}^4\bigr)^2
   +(A_{\rm sc}^2+4A_{\rm sc}^4)^2.
\end{equation}
Fix $N_0$ so large that $2N_0\ge b_p$ for every $p\in S_0$, and set
\begin{equation}\label{eq:Gamma}
 \Gamma=\prod_{p\in S_0}p^{N_0}.
\end{equation}
Then choose the fixed positive integer $D_0$ so that
\begin{equation}\label{eq:Dfixed}
 \prod_{p\in S_0}p\mid D_0\ \ \t{and}\ \ 
  D_0>\Gamma C_W/\sqrt\varepsilon.
\end{equation}

\begin{lemma}[Finite-local preparation]\label{lem:depth}
Fix a nonzero integer $r_0$ and put
$\mathcal T_{\rm loc}=S_0\cup\Supp(r_0)$.  There are a positive
integer $\kappa$, divisible by $c_*$ and by every prime in
$\mathcal T_{\rm loc}$, and sufficiently large integers
$a_q\ge v_q(\kappa)$ for $q\in\mathcal T_{\rm loc}$ such that:
\begin{enumerate}
\item if $h\equiv0\pmod\kappa$ and $v_q(h)=a_q$ for all
$q\in\mathcal T_{\rm loc}$, then $H=[h]P$ is sufficiently deep in
the formal group at every prime ideal $\mathcal Q\mid q$, and the two rational
components of $\varrho(H)/r_0$ belong to $q\mathbb Z_q$;
\item if $c\ne0$ is a sufficiently deep common multiple of $\kappa$
at these primes, with $v_q(c)-a_q$ sufficiently large, then for
$R_{\rm sc}=[c]P$ every normalized nonconstant leaf in
\eqref{eq:scalarleaves} belongs to $q\mathbb Z_q$, and
\begin{equation}\label{eq:Wunit}
 \mathcal W_{\rm sc}\equiv1\pmod q.
\end{equation}
\end{enumerate}
\end{lemma}
\begin{proof}
Choose $\kappa_q$ divisible by $q$ with $[\kappa_q]P$ in a
sufficiently deep formal subgroup at every prime ideal $\mathcal Q\mid q$; this is
possible by first entering the open formal subgroup and then
multiplying once more by $q$.  Take
$\kappa=\operatorname{lcm}(c_*,\{\kappa_q\})$, and only then
choose each $a_q\ge v_q(\kappa)$ sufficiently large.  Inside the
logarithm's isometric neighborhood,
\[
 v_{\mathcal Q}(z([h]P))=v_q(h)+\gamma_{q,\mathcal Q}
\]
with a fixed constant once $a_q$ is large enough, and the same formula
holds with $h$ replaced by $c$. Since $\varrho$ has a double zero at the identity, the first
assertion follows after increasing $a_q$; the same argument applied to
$[c]P$ controls the two coefficients of
$\varrho(R_{\rm sc})/r_0$.  Moreover
\[
 v_{\mathcal Q}\left(\frac{z([c]P)}{z([h]P)}\right)
 =v_q(c)-a_q.
\]
Taking this difference above $48v_q(r_0)+1$, and choosing the
formal neighborhoods coefficientwise inside
$q(\mathcal O_K\otimes_{\mathbb Z}\mathbb Z_q)$, makes both rational coefficients of the
ratio divided by $r_0^{48}$ divisible by $q$.  This coefficientwise
formulation is essential at the ramified prime $2$; membership only in
a local maximal ideal would not suffice.  The normalized leaves
therefore vanish modulo $q$, while the all-left constant leaf is one,
proving \eqref{eq:Wunit}.
\end{proof}

\subsection{The joint CRT lemma}

\begin{lemma}\label{lem:crt}
Let $\mathcal T_{\rm loc},\kappa,a_q$ be as above, with
$a_q\ge v_q(\kappa)$ and $a_2\ge6$. Fix the four root points and the
two scales $A_0,A_1$. Let $\mathcal S_{\rm sel}$ denote the finite set of
distinct rational primes produced from the contact positions after
repeated positions have been merged. The simultaneous-prime theorem
permits choosing $\mathcal S_{\rm sel}$ so that, with
\begin{equation}\label{eq:reductionmoduli}
 L_+=\lcm_{p\in\mathcal S_{\rm sel}}d_+(p)\ \ \t{and}\ \ 
 L_-=\lcm_{p\in\mathcal S_{\rm sel}}d_-(p),
\end{equation}
the system
\begin{align}
 h&\equiv0\pmod {L_+},&h&\equiv c_*\pmod {L_-},\label{eq:CRTa}\\
 h&\equiv0\pmod\kappa,&
 h&\equiv q^{a_q}\pmod{q^{a_q+1}}
 \quad(q\in\mathcal T_{\rm loc})\label{eq:CRTb}
\end{align}
has infinitely many integer solutions in one arithmetic progression.
A solution can be chosen arbitrarily large, with every ordinary
expression defined and bounded at both real embeddings, and with
$H=[h]P\in\mathcal O_H$.
\end{lemma}
\begin{proof}
For $0\le j\le1$ and $\theta\in\mathscr T$, put
$\alpha_{j,\theta}=-\theta/A_j$. Choose an integer $D_{\rm sh}>0$ such that
\[
 c_{j,\theta}=D_{\rm sh}^2\alpha_{j,\theta}\in\mathbb Z
 \qquad(0\le j\le1,\ \theta\in\mathscr T).
\]
These shifts are fixed because the four root points have already been
chosen. Remove repeated values before applying Sun's theorem, while
retaining their original occurrences in the three detectors.  The five
values with $j=0$ are distinct, so at least five distinct shifts remain
and the hypothesis $k\ge2$ is automatic.  Apply the theorem with
parameters $u,v$, and put $w=D_{\rm sh}v$ and $\eta=u/w^2$.
The selected prime occurrences are
\[
 p_{j,\theta}=u+c_{j,\theta}v^2
 =u-\frac{\theta}{A_j}w^2.
\]
The theorem permits prescribing a divisor of $v$, so every required
power of a prime in $S_0$ may be imposed on $w$. In the prime lower
bound exclude the finite supports of $D_{\rm sh}$, of the numerator and
denominator of every $A_j$, of all $\theta$ and $\theta-\xi$, of all
contact differences, of all fixed arithmetic exceptions, and of every
prime in $\mathcal T_{\rm loc}$. These exclusions depend only on the
already fixed roots and not on $h$.

Let $M_{\rm avoid}$ be the squarefree product of all odd prime divisors
of $\kappa r_0$, omitting $7$. Apply
\cite[Theorem 4.1]{Sun86} with this modulus, the distinct integral
shifts, sufficiently large prime lower bounds, and the prescribed
square-denominator dilation. It gives
\[
  \gcd(L_+,L_-)=4, \ v_2(L_+)=2,\ v_2(L_-)=6,\ 
  \gcd(L_+L_-,M_{\rm avoid})=1.
\]
The same theorem's Gaussian congruence at $7$ gives $7\nmid L_+L_-$.
Thus $L_+L_-$ is prime to the odd part of $\kappa$ and to every odd prime
in $\mathcal T_{\rm loc}$.

The first two congruences in \eqref{eq:CRTa} are compatible because
$4\mid c_*$. Their compatibility with $h\equiv0\pmod\kappa$ at odd
primes follows from the displayed coprimality, and at $2$ from
$64\mid c_*$. For $q\in\mathcal T_{\rm loc}$, the exact-valuation
congruence in \eqref{eq:CRTb} is compatible with
$h\equiv0\pmod\kappa$ because $a_q\ge v_q(\kappa)$. It is compatible
with the ordinary congruences at odd $q$ because $q\nmid L_+L_-$, and at
$2$ because $a_2\ge\max(6,v_2(\kappa))$ and $64\mid c_*$. The
generalized Chinese Remainder Theorem therefore gives one infinite
arithmetic progression of solutions. Lemma~\ref{lem:real} supplies
arbitrarily large positive terms of this progression in the required
real open set. Avoiding finitely many exceptional indices makes every
point and logarithm denominator nonzero.
\end{proof}

\begin{lemma}\label{lem:oddlucas}
Let $s,t,z\in\mathbb Z$, and suppose $a=|s|\ge7$, $z\ne0$, and
\begin{equation}\label{eq:lucastrue}
 t^2-s(s-2)tz+(s-2)^2z^2=1.
\end{equation}
Then we have
\begin{equation}\label{eq:lucasgrowthself}
 |t|>e^{a-3}.
\end{equation}
Conversely, if $s\ge7$ is odd and $s\equiv1\pmod3$, then there are
positive odd integers $t,z$ satisfying \eqref{eq:lucastrue}.
\end{lemma}
\begin{proof}
For the growth assertion, put
\[
 w=\operatorname{sgn}(s)(s-2)z,
 \qquad
 \rho=\frac{a+\sqrt{a^2-4}}2.
\]
Then \eqref{eq:lucastrue} becomes
\[
 t^2-atw+w^2=1.
\]
The element $t-\rho w$ is a norm-one unit of the order
$\mathbb Z[\rho]$.  After multiplying by $-1$ if necessary it is
positive in both real embeddings.  Every positive norm-one unit of
$\mathbb Z[\rho]$ is a power of $\rho$: indeed multiply such a unit by
a power of $\rho$ to obtain $v=A+B\rho$ with $1\le v<\rho$.  Since
$N(v)=1$, we have
\[
 B=\frac{v-v^{-1}}{\rho-\rho^{-1}}\in[0,1)\cap\mathbb Z,
\]
so $B=0$ and $v=1$.

Let $U_0=0$, $U_1=1$, and
$U_{k+1}=aU_k-U_{k-1}$ for $k=1,2,3,\ldots$.  Comparing the coefficients of
$\rho^{\pm k}$ gives
\[
 |w|=U_k\ \ \t{and}\ \ 
 |t|\in\{U_{k-1},U_{k+1}\}
\]
for some $k\ge1$.  Modulo $a-2$ one has $U_k\equiv k$, whereas modulo
$a+2$ one has $U_k\equiv(-1)^{k-1}k$.  Since
$|s-2|\mid w$, it follows that $k\ge a-2$ when $s>0$ and
$k\ge a+2$ when $s<0$.  The sequence is increasing for $a\ge3$, and
induction gives $U_j\ge(a-1)^{j-1}$ for $j\ge1$.  Therefore
\[
 |t|\ge(a-1)^{a-4}>e^{a-3}
\]
for $a\ge7$, proving \eqref{eq:lucasgrowthself}.

For existence, assume $s\ge7$ is odd and $s\equiv1\pmod3$.  Let
$U_0=0$, $U_1=1$, and $U_{j+1}=sU_j-U_{j-1}$ for $j=1,2,3,\ldots$. Set
$k=2(s-2)$,
\[
 t=U_{k+1},\qquad z=\frac{U_k}{s-2}.
\]
Modulo $s-2$ one has $U_j\equiv j$, so $z$ is an integer.  The known identity
\[
 U_{j+1}^2-sU_{j+1}U_j+U_j^2=1
\]
follows by induction and gives \eqref{eq:lucastrue}.  Modulo two the
recurrence has period $0,1,1$; since $k\equiv1\pmod3$, both $U_k$ and
$U_{k+1}$ are odd.  Finally $s-2$ is odd, so $z$ is odd, and the
recurrence is strictly increasing after $U_1$, giving $t,z>0$.
\end{proof}

\subsection{Joint completeness}

\begin{proposition}\label{prop:joint}
For every integer input $X$, there are nine points
\eqref{eq:ninepoints}, one rational slope $\eta$, two pairs of
ordinary pole parameters, one pair of scalar pole parameters, and all
finite odd-index and algebraic witnesses satisfying the complete
certificate. The same point $H$ occurs in all three pole certificates.
\end{proposition}
\begin{proof}
Fix $X\in\mathbb Z$, put $\beta=2X^2+3$, and take
$N=[\beta]P$. Then $\beta$ is a positive odd integer. Choose $\delta$
so that $r=\delta r_0(\beta)=|r_0(\beta)|>0$, and put
$\mathcal T_{\rm loc}=S_0\cup\Supp(r)$. Choose $\kappa,a_q$ as in
Lemma~\ref{lem:depth}, enlarging $a_2$ if necessary.

Choose a positive integer $c$, divisible by the fixed multiplier
$M$ of \eqref{eq:Mchoice}, by $6$, and by all prescribed local powers, and so deep at each
$q\in\mathcal T_{\rm loc}$ that Lemma~\ref{lem:depth}(2) will hold
whenever $v_q(h)=a_q$. The positive multiples of any fixed nonzero
integer have $G_{\mathbb R}^0$ in their orbit closure, so $c$ may also
be chosen arbitrarily large with
$R_{\rm sc}=[c]P\in\mathcal O_{\rm sc}$. Set
\[
 s=c+1,
 \qquad S=P+R_{\rm sc},
 \qquad Q=[c/M]P\ne O.
\]
Then $s\ge7$, $s$ is odd, and $s\equiv1\pmod3$. Apply
Lemma~\ref{lem:oddlucas} and put $T=[t]P$, $Z=[z]P$. Thus the four
root indices $\beta,s,t,z$ are positive and odd, and the odd-index
certificates of Lemma~\ref{lem:oddchart} have witnesses. The four
roots, and therefore the two scales, are now fixed.

Apply Lemma~\ref{lem:real} to the fixed roots. Invoke the simultaneous
prime theorem through Lemma~\ref{lem:crt}, choosing
$\eta=u/w^2$ sufficiently large and imposing the required small-prime
divisibility on $w$. This gives selected positive contact primes,
disjoint from $\mathcal T_{\rm loc}$, and a common $h$ satisfying
\eqref{eq:CRTa}--\eqref{eq:CRTb}. Put
\[
 H=[h]P,
 \qquad U_i=[hn_i]P\quad(1\le i\le4),
 \qquad(n_1,n_2,n_3,n_4)=(\beta,s,t,z).
\]
At each selected split contact prime, the reduced-group cardinalities put
$H,U_i$ in the formal kernel at one prime ideal $\mathcal P\mid p$ and
$V_0,V_i$ in the formal kernel at the conjugate prime ideal
$\mathcal P^\sigma$. By \cite[Proposition 2.3]{Burungale}, both
logarithm ratios agree with $n_i$ to precision $64$. Hence the first
fifteen tests have their required precision, and the final beta test
does as well because $n_1=\beta$.

We now verify the full hypotheses of pole completeness. Uniform real
boundedness is independent of the later progression modulus, and as
$\eta=u/w^2\to+\infty$ one has, for $0\le j\le1$,
\begin{align*}
 R_j&=\frac1{A_j^6\eta}(1+O(\eta^{-1})),
 &\mathcal W_j&=1+O(\eta^{-3}),\\
 \Theta_j&=2+O(\eta^{-3}),
 &B_j&=\frac{1+O(\eta^{-1})}{4A_j^6\eta}.
\end{align*}
For sufficiently large $\eta$, every $R_j$ is positive.  At each
real embedding every factor in $\mathcal W_j$ is then at least one,
so both conjugates of $\mathcal W_j$ and the trace $\Theta_j$ are
positive.  Thus both real nonzero guards hold and every $B_j$ is
positive and smaller than $\varepsilon$.  The exact factorization is
\begin{equation}\label{eq:factorR}
 R_j=\left(\frac1{w^5(t_j-\xi)^3}\right)^2
 \prod_{\theta\in\mathscr T}
 \left(u-\frac{\theta}{A_j}w^2\right).
\end{equation}
All five factors for a fixed detector are positive selected primes.
Suppose a selected prime $p$ divides $w=D_{\rm sh}v$. The lower bound excludes
the prime factors of $D_{\rm sh}$, so $p\mid v$. Since one selected factor
equals $p$, its defining equation gives $u\equiv0\pmod p$, and then
every selected value $u+c_{j,\theta}v^2$ is divisible by $p$. There
are at least two distinct shift values, because the five contact
positions belonging to $A_0=1$ are distinct. A positive prime divisible
by $p$ equals $p$, contradicting distinctness. Hence
\begin{equation}\label{eq:prime-coprime-w}
 \gcd\left(w,\prod_{j,\theta}p_{j,\theta}\right)=1.
\end{equation}

The lower bound also makes every selected prime a unit in the numerator
and denominator of every $A_j$, and in every contact difference. If
$p=p_{j,\theta}$, then $t_j\equiv\theta\pmod p$ and $t_j-\xi$ is a
unit. The five shifts for this fixed $j$ are distinct, so exactly one
factor of \eqref{eq:factorR} has valuation one and all other displayed
factors are units. Therefore
\begin{equation}\label{eq:Rj-valuation-one}
 v_p(R_j)=1\qquad(p=p_{j,\theta}).
\end{equation}
At every other prime outside the selected factors, $v_p(R_j)$ is even.
The norm--trace argument of \cite[Lemma 3.1]{Sun86} gives
$v_p(\Theta_j)=0$ at a selected prime: the successful split prime ideal
$\mathcal P\mid p$ contributes one to the trace residue, and its conjugate
prime ideal $\mathcal P^\sigma$ contributes zero or one, so the residue is
one or two. Thus $v_p(B_j)=1$ at each selected
prime for this detector, and every negative valuation outside $S_0$ is
even.  It remains to check the fixed exceptional primes; after that
check, these two facts together will imply $B_j\in\mathscr D_2$.

Fix $q\in S_0$, put
$\nu=v_q(w^2)$ and $\alpha_j=v_q(A_j)$, and make $\nu$ sufficiently
large. Since the zero-shift selected value is $u$, this $u$ is a prime
outside $S_0$, so $v_q(u)=0$. The valuation of
$t_j=A_ju/w^2$ is $\alpha_j-\nu$, smaller than the valuations of all
fixed contacts and of $\xi$. Hence
\[
 v_q(R_j)=\nu-6\alpha_j,
 \qquad v_q(t_j-\xi)=\alpha_j-\nu.
\]
All sixteen $D_i$ are integral at every prime ideal $\mathcal Q\mid q$ by
\cite[Lemma 8.1(1)]{Sun86}: the common index is divisible by $c_*$,
the root indices and $\beta$ are integers, and both anchors lie in the
prescribed formal neighborhoods. For a test of precision $k_i$,
\[
 v_{\mathcal Q}\left(
 R_j\left(\frac{\widetilde D_{j,i}}{R_j^{k_i}}\right)^2
 \right)
 \ge3\nu+(10k_i-8)\alpha_j
 \qquad(\mathcal Q\mid q).
\]
This tends to $+\infty$ with $\nu$, uniformly over the finite lists.
Thus $\mathcal W_j\equiv1$ to arbitrary $q$-adic precision at every
prime ideal $\mathcal Q\mid q$ and $\Theta_j\equiv2$ to sufficient
precision. Consequently
\[
 v_q(B_j)=\nu-6\alpha_j-2v_q(2)\ge b_q
\]
for all $q\in S_0$, after increasing the prescribed divisor of $w$.
It follows that each $B_j$ belongs to $\mathscr D_2$ and satisfies all
hypotheses of \eqref{eq:polecomplete}; that implication now supplies
both ordinary pole witnesses.

For the scalar certificate, Lemma~\ref{lem:depth} and the exact
condition $v_q(h)=a_q$ give
$\mathcal W_{\rm sc}\equiv1\pmod q$ for every
$q\in\mathcal T_{\rm loc}$, so this root is nonzero. If
$q\notin\Supp(r)$, then $v_q(B_{\rm sc})$ is even. If
$q\in\Supp(r)$, then $1+D_0r^2$ and $\mathcal W_{\rm sc}$ are local
units and
\[
 v_q(B_{\rm sc})=v_q(r)+2v_q(\Gamma)>0.
\]
Therefore $B_{\rm sc}\in\mathscr D_2$, and for $q\in S_0$ its
valuation is at least $2N_0\ge b_q$. The real bound
\eqref{eq:CW} and the choice \eqref{eq:Dfixed} give
\[
 0<B_{\rm sc}
 \le\frac{\Gamma^2C_W^2}{D_0^2r^3}
 \le\frac{\Gamma^2C_W^2}{D_0^2}<\varepsilon.
\]
Hence \eqref{eq:polecomplete} supplies the scalar pole witnesses.
The selected primes make the ordinary normalizers nonzero, the formal
conditions make both logarithm denominators nonzero, and the remaining
finite rational charts can be chosen on their regular loci. This proves
joint completeness.
\end{proof}

\section{The two height contradictions}\label{sec:heights}

Let
\begin{equation}\label{eq:Bbound}
 B=\max(1,|n|,|s|,|t|,|z|).
\end{equation}
For every nonzero index $m$ on the fixed curve,
\begin{equation}\label{eq:heightlower}
 \hW(x([m]P))=2\widehat h_E(P)m^2+O_E(1),
\end{equation}
with an implied constant depending only on the fixed curve and $P$.

\begin{lemma}[Ordinary height consequence]\label{lem:ordheight}
There is a fixed $C_{\rm ord}>0$ such that the two ordinary certificates imply
\begin{align}
 \Delta\ne0&\ \Longrightarrow\
 B^2\le C_{\rm ord}(1+|\Delta|^{3/8}),\label{eq:Deltabound}\\
 \beta\ne n&\ \Longrightarrow\
 B^2\le C_{\rm ord}(1+J^{3/2})
 \ \ \t{with}\ \ J=|\num(\beta-n)|.
 \label{eq:Jbound}
\end{align}
\end{lemma}
\begin{proof}
If $\Delta\ne0$, Lemma~\ref{lem:contact} and
\eqref{eq:finiteheight} give
\[
 \hW(A_j\eta)\le C_{\rm ht}\max(1,|\Delta|^{24/64})
 \qquad(0\le j\le1).
\]
Corollary~\ref{cor:twoscales} yields \eqref{eq:Deltabound}. If
$\beta\ne n$, the same argument uses the exponent $16$ and gives
$24/16=3/2$, proving \eqref{eq:Jbound}.
\end{proof}

\begin{lemma}[Numerator-height estimate]\label{lem:numerator}
There is a fixed constant $C_2>0$ such that every certificate with
$\beta\ne n$ satisfies
\begin{equation}\label{eq:logBbound}
 \log B\le C_2+3\hW(\beta).
\end{equation}
\end{lemma}
\begin{proof}
Write $\beta=a_0/b_0$ in lowest terms with $b_0>0$.  Since
$J=|a_0-nb_0|\ge1$, \eqref{eq:Jbound} gives, after enlarging a fixed
constant,
\[
 B^2\le C_3J^{3/2}.
\]
The definition of the rational height gives
$|a_0|,b_0\le e^{\hW(\beta)}$, and hence
\[
 J\le |a_0|+|n|b_0\le 2B e^{\hW(\beta)}.
\]
After enlarging the constant, substitution yields
$B^{1/2}\le C_4 e^{3\hW(\beta)/2}$ for a fixed $C_4>1$.
Taking logarithms proves the claim with $C_2=2\log C_4>0$.
\end{proof}

Choose the fixed positive integer $M$ so large that
\begin{equation}\label{eq:Mchoice}
\begin{gathered}
 6\mid M,\qquad M\ge8,\\
 \text{no real }B_0\ge M-1\text{ satisfies }
 B_0^2\le C_{\rm ord}(1+(10B_0^4)^{3/8}),\\
 C_2+3C_{\rm sc}\sqrt{c_0}<c_0-4
 \qquad\text{for every real }c_0\ge M,
\end{gathered}
\end{equation}
where $C_2$ is the constant of Lemma~\ref{lem:numerator}.
Such an $M$ exists because $B_0^2$ grows faster than $B_0^{3/2}$ and
$c_0$ grows faster than $\sqrt{c_0}$.

\begin{proposition}
\label{prop:certsound}
If the rational certificate has a solution at $X\in\mathbb Q$, then
$X\in\mathbb Z$.
\end{proposition}
\begin{proof}
Assume the rational certificate has a solution. The cyclic-point
encoding gives unique odd integers $n,s,t,z$. Write
$c=|s-1|$. Since $Q\ne O$ and $R_{\rm sc}=[M]Q$, one has
$c\ge M$ and $B\ge|s|\ge c-1\ge M-1$. Also $z\ne0$.

If $\Delta\ne0$, the elementary estimate $|\Delta|\le10B^4$ and
\eqref{eq:Deltabound} contradict \eqref{eq:Mchoice}. Hence
\begin{equation}\label{eq:Deltazero}
 \Delta=0.
\end{equation}
Lemma~\ref{lem:oddlucas} applies and gives
\begin{equation}\label{eq:Lucasgrowth}
 \log B\ge\log|t|>|s|-3.
\end{equation}

Suppose $\beta\ne n$. Lemma~\ref{lem:numerator},
Lemma~\ref{lem:scalarsound}, and \eqref{eq:Mchoice} give
\[
 \log B\le C_2+3C_{\rm sc}\sqrt c<c-4\le |s|-3,
\]
contrary to \eqref{eq:Lucasgrowth}. Thus $\beta=n\in\mathbb Z$.
If $X=a_X/b_X$ in lowest terms with $b_X>0$, then
$2a_X^2+3b_X^2=\beta b_X^2$, so $b_X^2\mid2a_X^2$.
Coprimality gives $b_X^2\mid2$, hence $b_X=1$ and
$X\in\mathbb Z$.
\end{proof}

\begin{proposition}
\label{prop:certcomplete}
If $X\in\mathbb Z$, then the rational certificate has a solution.
\end{proposition}
\begin{proof}
Proposition~\ref{prop:joint} supplies all points, continuous pole
parameters, finite odd-index roots, and every nonzero condition in the
rational certificate.
\end{proof}

\section{One final rational variable and the polynomial}\label{sec:compile}

We now justify the promised variable count and the passage to
an integer polynomial. The only continuously free rational quantities
are
\[
 \begin{array}{l|r}
 \text{two abscissa components for each of nine doubling preimages}&18\\
 \text{one common slope }\eta&1\\
 \text{two structural parameters for each of two ordinary }\Phi_0&4\\
 \text{two structural parameters for the scalar }\Phi_0&2
 \end{array}
\]
Their total is $25$.  Every remaining choice is finite algebraic over
these parameters and will be absorbed into one final rational unknown.

\subsection{Encoding the cyclic points}

A point of $2E(K)$ is the double of an affine point
$V\in E(K)$; even the identity has an affine two-torsion preimage.
Write $x(V)=u+v\rho$ with $u,v\in\mathbb Q$. The equation
$y(V)^2=x(V)^3-(11-6\rho)x(V)$ gives two rational polynomial
relations in the two rational components of the ordinate.
For fixed $(u,v)$, every allowed ordinate is a root of a fixed-degree
algebraic equation, and the distinct-root and zero-ordinate branches
can be separated by finitely many polynomial equalities and
inequalities. The explicit monic-separable splitting, including
the degenerate cases, is given in \cite[Section 10.1]{Sun86}.
Doubling, multiplication by $M$, addition, negation, conjugation,
and evaluation of the truncated logarithm are computed on a finite
collection of rational charts. Their equations and denominator guards
introduce no additional free coordinates. In particular
$R_{\rm sc}=[M]Q$ and $S=P+R_{\rm sc}$ are computed points rather
than additional freely represented points.

For each of $N,T,Z$, the odd-index certificate first chooses a
$K$-square root $h_Q$ of $x(Q)$, guarded by $h_Q\ne0$, then computes
$j_Q=y(Q)/h_Q$, and finally chooses the two $K$-square roots in
\eqref{eq:oddchart}.  An equation $U^2=A$ in $K$ has the same finite
four-branch rational description as Lemma~\ref{lem:ordinatecharts}.

\subsection{Complete finite chart catalogue and its verification}
\label{subsec:charts}
The assertion that nine cyclic points require only eighteen freely
chosen rational coordinates deserves an explicit algebraic proof.
The following are finite choice rules, not new quantified variables.

\begin{lemma}\label{lem:ordinatecharts}
Let $u,v\in\mathbb Q$, and set
\begin{align*}
 F_0&=u^3+6uv^2-11u+12v,\\
 F_1&=3u^2v+2v^3-11v+6u.
\end{align*}
The equation $(y_0+y_1\rho)^2=(u+v\rho)^3-(11-6\rho)(u+v\rho)$
is equivalent to exactly one of the following four chart types:
\begin{enumerate}
\item $F_1\ne0$: a rational root $a$ of
$a^4-F_0a^2+F_1^2/2=0$, with
$y_0=a$ and $y_1=F_1/(2a)$;
\item $F_1=0$, $F_0\ne0$, $y_1=0$: a rational root $a$ of
$a^2-F_0=0$, with $(y_0,y_1)=(a,0)$;
\item $F_1=0$, $F_0\ne0$, $y_0=0$: a rational root $a$ of
$a^2-F_0/2=0$, with $(y_0,y_1)=(0,a)$;
\item $F_0=F_1=0$: the unique root $a=0$ of $a=0$,
with $(y_0,y_1)=(0,0)$.
\end{enumerate}
Every displayed root polynomial is monic and separable on its stated
parameter branch; the four types exhaust all rational ordinate choices.
\end{lemma}
\begin{proof}
Comparison of rational components gives
$y_0^2+2y_1^2=F_0$ and $2y_0y_1=F_1$.
If $F_1\ne0$, then $y_0\ne0$ and elimination of $y_1$
gives the displayed quartic. Its derivative is
$2a(2a^2-F_0)$. A repeated root would force
$F_0^2=2F_1^2$, impossible for nonzero rational $F_1$.
If $F_1=0$ then $y_0y_1=0$; splitting according to which
coordinate vanishes gives the remaining three types. When
$F_0\ne0$ both quadratic polynomials are separable, while the
zero case is represented by the separable linear equation $a=0$.
\end{proof}

\begin{lemma}\label{lem:groupcharts}
For points $V,W\in E(K)$, the graph of their sum is covered by the
following finite rational charts: $V=O$; $W=O$; both affine with
distinct abscissas; both affine with equal abscissas and opposite
ordinates, giving $O$; and two equal affine points with nonzero
ordinate, using the tangent formula. The last chart is omitted when
the ordinate is zero, since such a point is two-torsion and doubles
to $O$. On the distinct-abscissa and tangent charts the corresponding
slopes are
\[
 \Lambda_{V,W}=\frac{y(W)-y(V)}{x(W)-x(V)},\qquad
 \Lambda_{V,V}=\frac{3x(V)^2-(11-6\rho)}{2y(V)}.
\]
Writing $\Lambda$ for the applicable slope on either of these
two charts, one has
$x(V+W)=\Lambda^2-x(V)-x(W)$ and
$y(V+W)=\Lambda(x(V)-x(V+W))-y(V)$.
Each denominator is nonzero on its chart. Consequently the graph
of any fixed integer multiplication $[m]$, in particular $[2]$,
$[M]$ and $[c_*]$, is given by a finite Boolean combination of
rational-function equations and nonzero-denominator guards.
\end{lemma}
\begin{proof}
If two affine points have equal abscissas, their ordinates satisfy
$(y(W)-y(V))(y(W)+y(V))=0$. These are precisely the inverse
and tangent cases, with the zero ordinate assigned to the inverse
case. The line-and-tangent formulas follow by substituting the line
into the Weierstrass equation and comparing its third intersection.
Direct substitution proves the output lies on the curve. Induction on
the standard binary double-and-add algorithm proves the assertion
for a fixed integer $m$, with a finite number of chart choices.
The identity and inverse cases are always retained, so the
exceptional substitutions at torsion points introduce no spurious
branches or missing points.
\end{proof}

For every one of the nine input abscissas we impose one of the four
ordinate charts, and then compute its double using the complete
group-law charts. All subsequent additions and fixed multiplications
are computed by the same rules. This proves that every resulting
cyclic point, including $O$, is covered, and conversely that every
choice accepted by a chart is a genuine point of $2E(K)$. In
particular, $O$ is attained by doubling the affine point $(0,0)$,
while the global guards $Q\ne O$ and $Z\ne O$ are imposed only on
the resulting doubled points. This distinction is essential.

The finite branch system further includes identity versus affine
branches for every occurrence of $\varrho$ or $z$; the three odd-index
certificates; the finite multiplier choice in each of the three copies
of $\Phi_0$; and, for each pole predicate's two square tests, either a
zero input or a nonzero input with a separable quadratic root equation.
The chart guards explicitly include
$Q,R_{\rm sc},S,N,T,Z,H,V_0\ne O$, the nonzero $h_Q$ in every odd chart, $\mathcal W_{\rm sc}\ne0$, the ordinary
$A_0A_1\ne0$, $R_j,\mathcal W_j,\Theta_j\ne0$, all
$t_j\notin\mathscr T\cup\{\xi\}$, both logarithm denominators, and
every denominator of a selected group-law or rational-function chart.
Field equalities and inequalities over $K$ are reduced to rational
components or field norms.
This is an exhaustive finite \emph{recursive specification}, including
specializations where a point becomes $O$ or an ordinate becomes zero;
its enormous expanded list need not be materialized to establish the
existential projection. The one-variable elimination in the next
subsection applies to this entire finite system, not to a selection
of generic charts.

\subsection{Finite root combination and zero branches}

The finite data consist of two different kinds.  The multiplier labels,
point-chart labels and zero/nonzero labels form a finite disjunction of
branches.  Within a fixed branch, the nine ordinates, the nine
quadratic-field square roots occurring in the three odd-index charts,
and the six rational pole square roots are roots of bounded-degree
polynomials over the rational function field in our $25$ free
parameters.  For nonzero rational square inputs one may adjoin the
quadratic roots separately.  Alternatively, the separable finite-root packaging of
\cite[Section 10.4]{Sun86}, based on Lemma~\ref{lem:ZS}, collects a
finite family of square conditions into one finite algebraic choice.
A square input equal to zero has its own finite branch and requires
no square-root variable. If no nontrivial algebraic root remains,
we include the harmless monic equation $T-1=0$.

Here are the details of the final compilation, following
\cite[Sections 10.3--10.6]{Sun86}.  On one branch, introduce the
ordinate and square-root choices successively.  After decomposing the preceding finite algebra into finitely many
components and branching on the corresponding idempotents, every new
coordinate satisfies a bounded-degree monic equation whose
discriminant is a unit on that stratum.  Independent choices may be
tensored, while dependent choices such as the two later roots in
\eqref{eq:oddchart} are adjoining steps in the same tower.  After
finitely many steps one obtains a finite \emph{étale} algebra over the
rational function field in the free parameters; here this means a finite
reduced algebra whose geometric factors are finite separable field
extensions.  A geometric point means a point after base change to a fixed
algebraic closure of the rational function field.

If the chosen branch has root coordinates
$r^{\rm alg}_1,\ldots,r^{\rm alg}_{N_{\rm alg}}$, choose first a
separating element
$\zeta_{\rm sep}=r^{\rm alg}_1+b_{\rm sep}r^{\rm alg}_2+\cdots
+b_{\rm sep}^{N_{\rm alg}-1}r^{\rm alg}_{N_{\rm alg}}$ for the full finite
\emph{étale} algebra,
using a finite predetermined list of integers $b_{\rm sep}$.  Write $Z_{\rm alg}$
for the characteristic-polynomial indeterminate.  The monic characteristic
polynomial $R_{b_{\rm sep}}(Z_{\rm alg})$ has nonzero discriminant, and the powers of
$\zeta_{\rm sep}$ express every coordinate $r_i^{\rm alg}$ as a rational
function of $Z_{\rm alg}$.  Translate the remaining equalities into
$F_1,\ldots,F_{N_{\rm eq}}\in\mathbb Q[Z_{\rm alg}]$ and the remaining
nonzero guards into $G_1,\ldots,G_{N_{\rm guard}}\in\mathbb Q[Z_{\rm alg}]$,
all reduced modulo $R_{b_{\rm sep}}$.  Taking monic gcds, we set
\[
 R_{\rm eq}=\gcd(R_{b_{\rm sep}},F_1,\ldots,F_{N_{\rm eq}})\ \ \t{and}\ \ 
 R_{\rm acc}=\frac{R_{\rm eq}}
 {\gcd(R_{\rm eq},\prod_{\ell=1}^{N_{\rm guard}}G_\ell)}.
\]
Because $R_{b_{\rm sep}}$ is separable, the roots of $R_{\rm acc}$ are exactly the
geometric tuples satisfying all equations and all nonzero guards.  A
constant $R_{\rm acc}$ rejects the branch.  Otherwise $R_{\rm acc}$ is again monic and
separable.  A rational root of $R_{\rm acc}$ corresponds to a unique geometric
tuple; Galois invariance and separation force every component of that
tuple to be rational.  Conversely, every accepted rational tuple gives
a rational root.  This is the ordering used in
\cite[Sections 10.3 and 10.5]{Sun86}; in particular, the gcd filtering
is performed only after the finite algebra has been put in a
separating one-variable form.

The degree, gcd, discriminant, and pivot conditions can all be
made uniform over the free parameters. Whenever a leading
coefficient might vanish, split into the finite branches recording
the actual degree, with the appropriate zero and nonzero tests.
Apply the same procedure to squarefree parts and to the
inverses in the finite algebra. The branch tree remains finite
because every initial polynomial has a fixed degree bound.
A branch with no remaining geometric point is discarded.
Let $\mathbf y\in\mathbb Q^{25}$ be the free-parameter tuple and let
$d_{\rm alg}$ be the degree of the retained monic polynomial on this
branch.  After choosing one common polynomial denominator
$D_{\rm den}(\mathbf y)$, write its polynomial coefficient numerators as
$\mathcal N_1(\mathbf y),\ldots,\mathcal N_{d_{\rm alg}}(\mathbf y)$.
Then the polynomial is
\begin{equation}\label{eq:monicpoly}
 \Pi_{\mathbf y}(T)=T^{d_{\rm alg}}+
 \sum_{j=1}^{d_{\rm alg}}
 \frac{\mathcal N_j(\mathbf y)}{D_{\rm den}(\mathbf y)}
 T^{d_{\rm alg}-j}.
\end{equation}
On the branch in question, a rational root is equivalent to a
rational tuple of the original finite choices satisfying all
constraints.  The separation discriminant is required to be nonzero;
without it a rational characteristic root need not represent a
rational tuple.  All branch coefficients and their denominators are
recorded as equations or inequations.  If a branch has no allowed
algebraic point, it is discarded.  Replacing $\Pi_{\mathbf y}(T)$ by
$\Pi_{\mathbf y}(T-c)$ for $c=0$ or $1$ gives two finite subbranches
and guarantees that at least one corresponding rational root is
nonzero whenever a rational root exists.

For completeness we recall the final reciprocal step.  Let
$G_{\rm guard}(\mathbf y)$ be the product of the numerators of the
remaining nonzero guards.  The condition
$D_{\rm den}(\mathbf y)\ne0$ is included among those guards.  Define
\begin{equation}\label{eq:reciprocal}
 \Psi(\mathbf y,W)
 =1+\sum_{j=1}^{d_{\rm alg}} \mathcal N_j(\mathbf y)
 D_{\rm den}(\mathbf y)^{2j-1}G_{\rm guard}(\mathbf y)^j W^j.
\end{equation}
If $D_{\rm den}(\mathbf y)G_{\rm guard}(\mathbf y)=0$, then
$\Psi(\mathbf y,W)=1$ for every $W$.  Otherwise, setting
$C_{\rm rec}(\mathbf y)=D_{\rm den}(\mathbf y)^2G_{\rm guard}(\mathbf y)$
gives the identity
\begin{equation}\label{eq:recipidentity}
 \Psi(\mathbf y,W)=
 (C_{\rm rec}(\mathbf y)W)^{d_{\rm alg}}
 \Pi_{\mathbf y}\!\left(\frac1{C_{\rm rec}(\mathbf y)W}\right).
\end{equation}
Moreover $\Psi(\mathbf y,0)=1$, so every solution has $W\ne0$.
Thus the final one-variable polynomial both encodes the nonzero
rational root and enforces every guard, including
$Q,R_{\rm sc},S,N,T,Z,H,V_0\ne O$,
$\mathcal W_{\rm sc}\ne0$,
the nonzero ordinary normalizers, and the two truncated-logarithm
denominators.  The finitely many equations in one branch may be combined by a
sum of squares, since a sum of rational squares is zero exactly when
every summand is zero.  The finite branches themselves may be combined
by a product.  The same $W$ is reused in each branch, since only
one branch needs to vanish.  Clearing fixed rational coefficients
produces one polynomial in $\mathbb Z[X,Y_1,\ldots,Y_{25},W]$;
relabeling $W$ as $Y_{26}$ gives the polynomial stated in
Theorem~\ref{thm:main}.

\begin{lemma}\label{lem:compiler}
Let $\mathbf a$ be any finite tuple of rational free parameters.
Suppose an existential condition on $\mathbf a$ consists of finitely
many branches, each with finitely many equations and inequations in
rational functions of $\mathbf a$ and in algebraic auxiliary roots
introduced by a finite tower of bounded-degree monic equations.
Suppose each adjoining equation can be split into monic separable
equations on finitely many conditions on the preceding finite algebra.
Then the whole condition is expressible by one
integer-coefficient polynomial equation in $\mathbf a$ and one
further rational unknown, without changing its projection.
\end{lemma}
\begin{proof}
Work over one rational specialization of the free parameters.  On a
fixed stratum, successively adjoining the roots gives a finite
\emph{étale} algebra, meaning here a finite reduced algebra whose
geometric factors are finite separable field extensions.  The unit-discriminant
condition is imposed componentwise, not merely as a nonzero element of a
product algebra.  Let
$D_{\rm alg}$ be the rank of this finite étale algebra, equivalently
the number of its geometric points, and let $r_{\rm alg}$ be the number
of auxiliary coordinates on the stratum.  Both admit uniform bounds
from the fixed degree tower.  For a coefficient vector
$\mathbf c=(c_1,\ldots,c_{r_{\rm alg}})$, the linear combination
$\sum c_i\alpha_i$ fails to separate a fixed pair of distinct geometric
points only when $\mathbf c$ lies on one nonzero linear hyperplane.
Choose an integer $L$ larger than the uniform bound for the number of
pairs.  In the grid $\{0,\ldots,L\}^{r_{\rm alg}}$, each such
hyperplane contains at most $(L+1)^{r_{\rm alg}-1}$ points, whereas the
whole grid has $(L+1)^{r_{\rm alg}}$ points.  The union bound therefore
leaves a coefficient vector separating all geometric points.  Thus a
fixed finite grid of integer linear combinations contains a separating
choice uniformly on every stratum.  Express all auxiliary coordinates
in the resulting one-variable algebra.  Polynomial equalities and inequations then
select a Galois-stable subset by the univariate gcd construction
written above.  If the subset is empty, reject the branch.  Otherwise
the retained monic squarefree factor has rational roots exactly when
the selected set has a rational point: a rational separating value is
fixed under Galois and uniqueness forces the full tuple to be fixed.
Conversely any rational tuple yields a rational root.

To obtain this construction uniformly in $\mathbf a$, branch on the
actual degrees of all specialized polynomials, their gcds and
squarefree factors, the leading coefficients, and the pivots in the
finite algebra. Every such condition is an equality or inequation in
rational functions of $\mathbf a$, and the number of branches is
finite because all initial degrees are bounded. Before the
reciprocal conversion, replace the retained monic polynomial by $\Pi(T-c)$ for
$c=0$ or $1$ on two finite subbranches, so that one accepted root is
nonzero whenever a rational root exists. Record the nonzero guards and
use \eqref{eq:reciprocal} to convert that nonzero rational root into a
single polynomial equation.
Finally combine complete branch conjunctions by sums of squares
and their finite disjunction by products. The same final variable
is reused in every branch. This is the explicit uniform
specialization procedure used in \cite[Sections 10.3--10.6]{Sun86}.
\end{proof}

\begin{remark}
Once the fixed arithmetic data and one admissible integer $M$ are
supplied, the finite-algebra prescription is a finite algorithm, although
it may be enormous.
Its existence and its projection equivalence, rather than a
monomial-expanded coefficient table or a degree estimate, are what
Theorem~\ref{thm:main} asserts.  Lemmas~\ref{lem:oddchart},
\ref{lem:ordinatecharts}, \ref{lem:groupcharts}, and
\ref{lem:compiler}
give an exhaustive
mathematical branch specification. They do not amount to an
expanded coefficient file or a machine-formal verification of all
compiled specializations; no such computational certificate is
claimed here. The shared slope and elliptic
indices are not additional root choices; the latter appear only
in the arithmetic proof.
\end{remark}

\section{A compact exact certificate for the fixed multiplier}
\label{sec:multdag}

The fixed multiplier $M$ from \eqref{eq:Mchoice} may be extraordinarily
large.  It is important to distinguish three different mathematical objects:
(i) the integer $M$ fixed by the height inequalities; (ii) a rational-function
arithmetic circuit computing $[M]$; and (iii) the fully expanded numerator,
denominator, and final coefficient table.  Only (i) and (ii) are needed for
the existential construction.  In particular, a verified circuit identity
can serve as an exact finite specification without producing all monomials.

\begin{lemma}\label{lem:multdag}
Let $L$ be a field of characteristic zero, $A\in L$, and
$E_A:y^2=x^3+Ax$.  For an affine point $R=(x,y)$ with $y\ne0$,
define $\psi_0=0$, $\psi_1=1$, $\psi_2=2y$, and
\begin{align*}
 \psi_3&=3x^4+6Ax^2-A^2,\\
 \psi_4&=4y(x^6+5Ax^4-5A^2x^2-A^3).
\end{align*}
For every $k\ge2$, and for $k\ge3$ in the even-index formula,
use the recurrences
\begin{align}
 \psi_{2k+1}&=\psi_{k+2}\psi_k^3-\psi_{k-1}\psi_{k+1}^3,\label{eq:psi_odd}\\
 \psi_{2k}&=\frac{\psi_k}{2y}
       (\psi_{k+2}\psi_{k-1}^2-\psi_{k-2}\psi_{k+1}^2).
       \label{eq:psi_even}
\end{align}
Here the recurrences are evaluated in the function field of $E_A$;
all indicated quotients by $2y$ have polynomial meaning after the
curve equation is used.  Put
\[
 \phi_m=x\psi_m^2-\psi_{m+1}\psi_{m-1}.
\]
For every integer $m\ge1$ with $[m]R\ne O$, we have
\begin{equation}\label{eq:multxy}
 x([m]R)=\frac{\phi_m(R)}{\psi_m(R)^2}\ \ \t{and}\ \ 
 y([m]R)=\frac{\psi_{2m}(R)}{2\psi_m(R)^4}.
\end{equation}
The recurrences give a finite directed acyclic arithmetic circuit for
these expressions for every fixed $m$.  On the domain $R\in\mathbb ZP\setminus\{O\}$ of
Theorem~\ref{thm:main}, none of the denominators in
\eqref{eq:multxy} vanishes.
\end{lemma}
\begin{proof}
The recurrences, the formula $x([m]R)=\phi_m/\psi_m^2$, and the standard identity $\psi_{2m}=2\psi_m\omega_m$, where
$\omega_m$ is the usual ordinate numerator division polynomial, with
$y([m]R)=\omega_m/\psi_m^3$ are the classical elliptic division-polynomial identities in
\cite[Section 2.2]{Glazer}.  They give
\eqref{eq:multxy} by substitution.  They also identify zeros of
$\psi_m$ with points satisfying $[m]R=O$.  Since the chosen
subgroup $\mathbb ZP$ is infinite cyclic and torsion-free,
$R\ne O$ implies $[m]R\ne O$ and $y(R)\ne0$ for every
$m\ge1$.  Thus every division used above is valid at the actual
characteristic-zero points.

Memoizing every division polynomial that occurs in the terminating
recursion gives a finite arithmetic circuit for the fixed integer
$m$.  No asymptotic circuit-size estimate is needed for the variable
count or for the projection argument.
\end{proof}

\begin{remark}
\label{rem:mult_scope}
The choice \eqref{eq:Mchoice} specifies $M$ by
inequalities involving fixed height-comparison constants. 
 The circuit of Lemma~\ref{lem:multdag} is an
exact algebraic specification for every such choice of $M$.
  A program instantiated at a
sample value of $M$ could test the division-polynomial implementation
and exceptional modular reductions, but such tests would not constitute
an independent formal verification of the full twenty-six-unknown
projection.  Any separately distributed implementation should therefore
state its chosen $M$ and carry its own version and hash information.
\end{remark}

\section{Explicit control of the rational projection}
\label{sec:projection}

For clarity we isolate the projection assertion implicit in
Lemma~\ref{lem:compiler}.  The following result is the precise
logical bridge from the finite atlas to a polynomial with twenty-six
existential unknowns. Its strength is the \emph{all-specializations}
qualification; generic separability alone would be insufficient.

\begin{proposition}
\label{prop:exactprojection}
Fix all arithmetic constants, including one integer $M$ satisfying
\eqref{eq:Mchoice}. Let $\mathcal R_M(X;\mathbf y;\mathbf t)$ be the
finite Boolean combination of rational equations and inequations assembled
in Sections~\ref{sec:ordinary}--\ref{sec:scalar} together with the
point, chart and finite-algebra conditions of Section~\ref{sec:compile},
where $\mathbf y$
consists of exactly $25$ free rational parameters and $\mathbf t$
denotes the finite algebraic ordinate, odd-chart, and pole-root choices.
Each rational component of $\mathbf t$ is required to be rational in an
accepted solution.  We may take $N_{\rm aux}=24$: nine ordinate roots,
nine quadratic-field square-root coordinates from the three odd-index
charts, and six rational pole-square roots.  There is an integer polynomial
\[
 F_M\in\mathbb Z[X,Y_1,\ldots,Y_{25},W]
\]
for which, at every rational value $x$,
\begin{equation}\label{eq:fullprojection}
 \begin{split}
 &\exists\,\mathbf y\in\mathbb Q^{25}
   \;\exists\,\mathbf t\in\mathbb Q^{N_{\rm aux}}
   \,[\mathcal R_M(x;\mathbf y;\mathbf t)]\\
 &\hspace{15mm}\Longleftrightarrow
   \exists\,\mathbf y\in\mathbb Q^{25}
   \;\exists\,w\in\mathbb Q\,
   [F_M(x,\mathbf y,w)=0].
 \end{split}
\end{equation}
Branches using fewer auxiliary coordinates are padded by equations
$t_i=0$.  These coordinates belong to the finite algebra encoded by
$W$ and are not additional free parameters in the final polynomial.
Note that 
\eqref{eq:fullprojection} includes every exceptional zero,
denominator, and non-separable specialization.
\end{proposition}
\begin{proof}
We spell out the input and the filtering steps of the finite-
separable-algebra algorithm of \cite[Sections 10.1--10.6]{Sun86}.
First branch on the ordinate and group-law charts, the three odd-index
charts, one finite multiplier label for each of the three instances of $\Phi_0$, and
the zero/nonzero cases of the six rational pole-square inputs.  Each
odd chart uses three $K$-square equations: one for $h_Q$, one for $u_Q$,
and one for $v_Q$; after $h_Q\ne0$ has been imposed, $j_Q=y(Q)/h_Q$
is a guarded rational function.  By the same calculation as
Lemma~\ref{lem:ordinatecharts}, each $K$-square equation is represented
by one monic equation of degree at most four on finitely many branches.
Thus, if the six pole square conditions are left uncombined for
this audit, a branch has at most eighteen quartic root coordinates and
six quadratic root coordinates, and its geometric root count is at
most $4^{18}2^6=2^{42}$.  The Zhang--Sun packaging \cite[Theorem 1.2]{ZS} gives an equivalent
smaller presentation, but is not needed for this upper bound.  This is
only a finite-algebra bound, not a bound on the degree of $F_M$.

All field-valued expressions in $\mathbb Q(\sqrt2)$ are
split into their two rational components.  An equality in the field
is equivalent to two rational equalities; a nonvanishing condition
is equivalent to the nonvanishing of its rational field norm.
Every rational-function denominator is retained as an explicit
nonzero guard, including those in the fixed-multiplication circuit
and the logarithm ratios.  By Lemma~\ref{lem:multdag}, the new
$[M]Q$ circuit is defined on every accepted branch with $Q\ne O$.
The adjoining equations are made monic and separable after splitting
on their actual degree and discriminant, including every zero
specialization. There are finitely many branches because the tower
length and all degrees are bounded independently of the free
parameters.

In a fixed branch and at a fixed rational specialization of
$(x,\mathbf y)$, form the resulting finite separable tower, tensoring
independent blocks. The result is finite \emph{étale}, hence reduced, with one
geometric point for each tuple of root choices.  Let
$r^{\rm alg}_1,\ldots,r^{\rm alg}_{N_{\rm br}}$ be its chosen root
coordinates.  Before
filtering by the remaining equations and guards, choose a separating
linear combination
\[
 \zeta_{\rm sep}=r^{\rm alg}_1+b_{\rm sep}r^{\rm alg}_2+\cdots
 +b_{\rm sep}^{N_{\rm br}-1}r^{\rm alg}_{N_{\rm br}}
\]
from a fixed finite list of integers $b_{\rm sep}$.  Two distinct geometric
points collide for at most $N_{\rm br}-1$ values of $b_{\rm sep}$; hence one
choice has a characteristic polynomial $R_{b_{\rm sep}}(Z_{\rm alg})$ with nonzero
discriminant.  On that branch, the Vandermonde recovery described in
\cite[Section 10.3]{Sun86} expresses every $r_i^{\rm alg}$ as a rational function of
$Z_{\rm alg}$ modulo $R_{b_{\rm sep}}$.

Translate every remaining equality and nonzero guard into a univariate
polynomial in $Z_{\rm alg}$.  Successive squarefree gcds with the equality
polynomials, followed by removal of the factors on which the product
of the guards vanishes, produce a monic squarefree factor $R_{\rm acc}(Z_{\rm alg})$.
Its roots are exactly the geometric tuples accepted by the branch.
The degree, gcd, discriminant and pivot operations are made uniform by
finitely branching on their actual values, as in
\cite[Sections 10.5--10.6]{Sun86}.  The retained root set is
Galois-stable.  A rational root of $R_{\rm acc}$ is fixed by Galois, and
separation forces the corresponding whole tuple to be rational.
Conversely, any accepted rational tuple determines a rational root.
Replacing $R_{\rm acc}(Z_{\rm alg})$ by $R_{\rm acc}(Z_{\rm alg}-c)$ for $c=0$ or $1$ ensures that the
selected rational root is nonzero in at least one of the two finite
branches.

Finally write the retained monic polynomial as
\[
 \Pi(T)=T^{d_{\rm alg}}+\sum_{j=1}^{d_{\rm alg}}
       \frac{\mathcal N_j}{D_{\rm den}}T^{d_{\rm alg}-j},
\]
where $d_{\rm alg}$ is its degree, the $\mathcal N_j$ are polynomial
numerators, and $D_{\rm den}$ is one common polynomial denominator.
Let $G_{\rm guard}$ be the product of the numerators of every remaining
nonzero guard, including the guard $D_{\rm den}\ne0$.  Use
\[
 1+\sum_{j=1}^{d_{\rm alg}}\mathcal N_j D_{\rm den}^{2j-1}
       G_{\rm guard}^j W^j
 =\bigl(D_{\rm den}^2G_{\rm guard}W\bigr)^{d_{\rm alg}}
  \Pi\!\left(\bigl(D_{\rm den}^2G_{\rm guard}W\bigr)^{-1}\right)
\]
when $D_{\rm den}G_{\rm guard}W\ne0$.  If
$D_{\rm den}G_{\rm guard}=0$, the left side is identically $1$; its
constant term also shows that every solution has $W\ne0$.  This is
an exact polynomial equivalence for each specialization, not merely for
generic parameters.  Combine the
finite conjunctions by sums of squares and disjunctions by products;
the same variable $W$ is reused for every branch.  Clearing fixed
rational coefficients gives $F_M\in\mathbb Z[X,Y_1,\ldots,Y_{25},W]$.
The construction has not quantified any additional rational
parameter, proving \eqref{eq:fullprojection}.
\end{proof}

\begin{proof}[Proof of Theorem~\ref{thm:main}]
Choose and fix one integer $M$ satisfying \eqref{eq:Mchoice}.
Propositions~\ref{prop:certsound} and~\ref{prop:certcomplete} show that
the rational-input projection of the full certificate is exactly
$\mathbb Z$.  Proposition~\ref{prop:exactprojection} supplies an
integer polynomial $F_M$ in the input and $26$ existential rational
unknowns with exactly the same projection.  Taking $F=F_M$ proves the
theorem.
\end{proof}

\begin{remark}
\label{rem:projection_scope}
Proposition~\ref{prop:exactprojection} proves the projection
identity by a uniform algebraic procedure, assuming the source
finite-separable-algebra theorem as stated.  
\end{remark}

\section{An \texorpdfstring{$abc$}{abc}-conditional sixteen-unknown refinement}\label{sec:abc}

We prove Theorem~\ref{thm:abc16}. The point of the $abc$ conjecture
is that it replaces the power-law five-contact bound
\eqref{eq:fiveheight} by a logarithmic bound. Once this is available,
the Lucas amplification used in Section~\ref{sec:heights} is no longer
needed.

\begin{lemma}[Granville's $abc$ consequence]\label{lem:granville}
Assume the $abc$ conjecture. Let
$G(U,V)\in\mathbb Z[U,V]$ be a homogeneous binary form of degree
$r\ge3$ with no repeated factor. For every $\varepsilon>0$,
\[
 \max\{|u|,|v|\}^{r-2-\varepsilon}
 \ll_{G,\varepsilon}\operatorname{rad}G(u,v)
\]
for coprime integers $u,v$ with $G(u,v)\ne0$.
\end{lemma}
\begin{remark}
This is a consequence of the $abc$-conjecture proved by Granville
\cite{Granville}; see also the use of the same estimate in
\cite[Section 5]{KPV}.
\end{remark}

\begin{lemma}\label{lem:abcfive}
Assume the $abc$-conjecture, and fix five distinct rational points
$\mathscr A\subset\mathbb P^1(\mathbb Q)$ and a finite set $S$ of
primes. For $a=(u:v)\notin\mathscr A$, represented by coprime
integers, let
\[
 M_{\mathscr A,S}(a)=
 \prod_{\substack{p\notin S\\
 d_{p,c}(a)>0\text{ is odd for some }c\in\mathscr A}}p,
\]
where the normalized depths are defined as in \eqref{eq:contactdepth}
using primitive integral contact forms. Then there is a fixed constant
$C_{\mathscr A,S}$ such that
\begin{equation}\label{eq:abcfive}
 \hW(a)\le 2\log\max\{1,M_{\mathscr A,S}(a)\}
 +C_{\mathscr A,S}.
\end{equation}
\end{lemma}
\begin{proof}
Let $L_c(U,V)$ be the five primitive contact forms and put
$G(U,V)=\prod_{c\in\mathscr A}L_c(U,V)$. Enlarge $S$ by the fixed
primes dividing the pairwise resultants and the contents of the chosen
forms. Then for $p\notin S$ at most one $L_c(u,v)$ is divisible by
$p$. If $e=v_p(G(u,v))>0$, the contribution $p^2$ to
$\operatorname{rad}(G(u,v))^2$ is bounded by the $p$-part of
$|G(u,v)|M_{\mathscr A,S}(a)$: for even $e$ this follows from
$e\ge2$, and for odd $e$ from the extra factor $p$ in
$M_{\mathscr A,S}(a)$. The fixed primes contribute only a fixed
constant. Hence, with $H=\max\{|u|,|v|\}$,
\[
 \operatorname{rad}(G(u,v))^2
 \ll_{\mathscr A,S}|G(u,v)|M_{\mathscr A,S}(a)
 \ll_{\mathscr A,S}H^5M_{\mathscr A,S}(a).
\]
Lemma~\ref{lem:granville} with $r=5$ and $\varepsilon=1/4$ gives
$H^{11/4}\ll\operatorname{rad}(G(u,v))$. Squaring and comparing
with the preceding display yields $H^{1/2}\ll M_{\mathscr A,S}(a)$.
Taking logarithms proves \eqref{eq:abcfive}, after enlarging the
constant to include the case $M_{\mathscr A,S}(a)=1$.
\end{proof}

We now define the smaller certificate. Fix a positive integer
$M_{abc}$, divisible by $6$, to be chosen after the constants below are
known. Keep only four freely represented cyclic points
\begin{equation}\label{eq:abcpoints}
 Q,N,H,U\in\mathbb ZP.
\end{equation}
Put
\[
 \beta=2X^2+3,\qquad
 R_{\rm sc}=[M_{abc}]Q,\qquad S=P+R_{\rm sc},
\]
and write $N=[n]P$, $S=[s]P$. Require
$Q,R_{\rm sc},S,N,H\ne O$. Set
\[
 V_0=H-[c_*]P,\qquad V=U-[c_*]N,
\]
require $V_0\ne O$, and define
\[
 L=\frac{\ell_{64}(z(U))}{\ell_{64}(z(H))},\qquad
 L'=\frac{\ell_{64}^{\sigma}(z(V^\sigma))}
          {\ell_{64}^{\sigma}(z(V_0^\sigma))}.
\]
The two denominators are required to be nonzero. The ordinary test
list is
\begin{equation}\label{eq:abctests}
 \varrho(H),\quad\varrho(U),\quad
 \varrho(V_0^\sigma),\quad\varrho(V^\sigma),\quad
 L-L',\quad L-\beta,
\end{equation}
with required precisions
\begin{equation}\label{eq:abcorders}
 1,1,1,1,64,16.
\end{equation}
No Lucas points and no odd-index certificates occur in this smaller
system.

For the two nonzero root indices define
\begin{equation}\label{eq:abcaggregate}
 T_{abc}=\frac12\Tr_{K/\mathbb Q}(x(N)+x(S)),\qquad
 A_0=1,\quad A_1=T_{abc}^2.
\end{equation}
The sum form of the aggregate trace-height theorem
\cite[Section A38]{Glazer}, now with $r=2$, gives fixed constants
$C_{ag}>0$ and $C_{ag,0}\ge0$ such that
\begin{equation}\label{eq:abcaggregateheight}
 B^2\le C_{ag}\hW(T_{abc})+C_{ag,0},\qquad
 B=\max\{1,|n|,|s|\}.
\end{equation}
As in Lemma~\ref{lem:aggregatetrace}, $T_{abc}>0$. Introduce one
free slope $\eta$ and put $t_j=A_j\eta$. Define $R_j$ by
\eqref{eq:Rj}, and define the two products $\mathcal W_j$ exactly as
in \eqref{eq:Wj}, but with the six tests
\eqref{eq:abctests} and the six orders \eqref{eq:abcorders}. Define
$\Theta_j,B_j$ by \eqref{eq:Bj}, impose the same nonzero guards, and
require $\Phi_0(B_j)$ for $j=0,1$. Use exactly the scalar certificate
of Section~\ref{sec:scalar}, with $M$ replaced by $M_{abc}$.

\begin{lemma}\label{lem:abcordinary}
For every prime factor $p$ of $\mathfrak M(t_j)$, there is a
prime ideal $\mathcal P\mid p$ satisfying all six conditions
\eqref{eq:abctests} to the precisions \eqref{eq:abcorders}. Hence,
if $\beta\ne n$ and $J=|\num(\beta-n)|$, then
\begin{equation}\label{eq:abcJdiv}
 \mathfrak M(t_j)^{16}\mid J\qquad(j=0,1).
\end{equation}
\end{lemma}
\begin{proof}
The norm--trace argument in Lemma~\ref{lem:contact} applies verbatim
to any finite test list: if every place above $p$ failed at least one
test, the inverse factors would force $v_p(\Theta_j)\ge v_p(R_j)$,
contradicting the negative-odd-pole obstruction \eqref{eq:polesound}.
Thus one place satisfies all six tests simultaneously. The first four
make $H,U,V_0^\sigma,V^\sigma$ formal. Proposition 2.3 of
\cite{Burungale}, in the two-anchor calculation of
Lemma~\ref{lem:tw_anchor}, gives
\[
 L\equiv n\pmod{\mathcal P^{64}}.
\]
The last test therefore gives $v_p(\beta-n)\ge16$. Positive valuation
forces the reduced denominator of $\beta-n$ to be a $p$-unit, so
$p^{16}\mid J$. Multiplication over the distinct prime factors proves
\eqref{eq:abcJdiv}.
\end{proof}

\begin{proposition}\label{prop:abcsound}
Assume the $abc$-conjecture. If the smaller certificate above has a
rational solution, then $X\in\mathbb Z$, provided that $M_{abc}$ is chosen
sufficiently large in terms of the fixed arithmetic data.
\end{proposition}
\begin{proof}
Suppose first that $\beta\ne n$ and put
$J=|\num(\beta-n)|\ge1$. By Lemmas~\ref{lem:abcfive} and
\ref{lem:abcordinary}, for $j=0,1$ we have
\begin{equation}\label{eq:abctheight}
 \hW(t_j)\le\frac18\log J+C_0.
\end{equation}
Since $t_1/t_0=T_{abc}^2$, the usual height inequality gives
\[
 \hW(T_{abc})\le\frac12(\hW(t_0)+\hW(t_1))
 \le\frac18\log J+C_0.
\]
Thus \eqref{eq:abcaggregateheight} gives
\begin{equation}\label{eq:abcBfromJ}
 B^2\le C_1+\frac{C_{ag}}8\log J.
\end{equation}

Write $c=|s-1|$. Lemma~\ref{lem:scalarsound} still gives
$M_{5,S_{\rm sc}}(\beta)^{48}\mid c$. Applying
Lemma~\ref{lem:abcfive} to the five contacts $\mathscr C$ yields
\begin{equation}\label{eq:abcbetaheight}
 \hW(\beta)\le C_2+\frac1{24}\log c.
\end{equation}
If $\beta=a/b$ in lowest terms, then
\[
 J=|a-nb|\le(B+1)e^{\hW(\beta)},\qquad
 c=|s-1|\le B+1.
\]
Consequently,
\[
 \log J\le C_2+\frac{25}{24}\log(B+1).
\]
Substitution in \eqref{eq:abcBfromJ} gives
\begin{equation}\label{eq:abclogcontradiction}
 B^2\le C_3+\frac{25C_{ag}}{192}\log(B+1).
\end{equation}
All constants here depend only on the fixed curve, contact sets and
pole predicates; they are fixed before $M_{abc}$ is chosen. Since
$Q\ne O$, if $Q=[m_Q]P$ then $s-1=M_{abc}m_Q$ and therefore
$B\ge|s|\ge M_{abc}-1$. Choose $M_{abc}$, divisible by $6$, so large
that \eqref{eq:abclogcontradiction} has no solution with
$B\ge M_{abc}-1$. This contradiction proves $\beta=n\in\mathbb Z$.
Finally, if $X=a_X/b_X$ is in lowest terms with $b_X>0$, then
$2a_X^2+3b_X^2=\beta b_X^2$ implies $b_X^2\mid2a_X^2$, hence
$b_X=1$. Thus $X\in\mathbb Z$.
\end{proof}

\begin{proposition}\label{prop:abccomplete}
Under the $abc$-conjecture, each $X\in\mathbb Z$ has a solution of
the smaller certificate.
\end{proposition}
\begin{proof}
Put $\beta=2X^2+3$ and $N=[\beta]P$. Repeat the scalar preparation of
Proposition~\ref{prop:joint}, now choosing a positive integer $c$
divisible by $M_{abc}$, by $6$, and by the required fixed local powers,
with $R_{\rm sc}=[c]P$ in the same fixed real neighborhood. Set
\[
 Q=[c/M_{abc}]P,\qquad S=P+R_{\rm sc}=[c+1]P.
\]
The two root points $N,S$ and hence the two scales
\eqref{eq:abcaggregate} are now fixed. Apply the proof of
Lemma~\ref{lem:crt} to these two scales after deleting the unused
roots and tests. The simultaneous-prime theorem is stated for any
fixed finite set of shifts, and the five shifts belonging to $A_0=1$
are already distinct; hence no new prime-pattern assertion is needed.
It supplies a slope $\eta=u/w^2$ and one common comparison index $h$.
Put
\[
 H=[h]P,\qquad U=[h\beta]P.
\]
The same CRT congruences put $H,U$ in the formal kernel at the selected
prime ideal and $V_0,V$ in the formal kernel at its conjugate. By
\cite[Proposition 2.3]{Burungale}, the two logarithm ratios both agree
with $\beta$ to precision $64$, so all six tests are satisfied.

The proof of Proposition~\ref{prop:joint} now applies with fewer test
factors. In particular, the exact factorization \eqref{eq:factorR},
the selected-prime unit calculation, and the fixed-prime depth estimate
show that each ordinary input $B_j$ lies in $\mathscr D_2$, is positive
and smaller than $\varepsilon$, and has the required depths at $S_0$.
Hence \eqref{eq:polecomplete} supplies the two ordinary pole witnesses.
Lemma~\ref{lem:depth} and the unchanged scalar argument supply the
scalar pole witnesses. Removing tests and root points only removes
conditions from the finite preparation; it creates no new denominator
or congruence requirement. Thus all guards and finite chart choices can
be satisfied.
\end{proof}

\begin{proof}[Proof of Theorem~\ref{thm:abc16}]
Propositions~\ref{prop:abcsound} and \ref{prop:abccomplete} show that
the smaller rational certificate projects exactly to $\mathbb Z$.
Its continuously free rational parameters are
\[
 \underbrace{2\cdot4}_{Q,N,H,U}
 +\underbrace{1}_{\eta}
 +\underbrace{4}_{\text{two ordinary }\Phi_0}
 +\underbrace{2}_{\text{scalar }\Phi_0}=15.
\]
There are no Lucas points and no odd-index certificates. The remaining
ordinate choices, the six pole square roots, finite multiplier labels,
point-chart choices and nonzero guards form the same kind of finite
separable algebraic system treated by Lemma~\ref{lem:compiler}; with
four rather than nine freely represented cyclic points the construction
is strictly smaller. Lemma~\ref{lem:compiler} therefore uses one
further rational unknown and produces one integer-coefficient
polynomial with exactly $15+1=16$ existential unknowns and projection
$\mathbb Z$.  This concludes the proof.
\end{proof}

\end{document}